# NORMAL APPROXIMATION UNDER LOCAL DEPENDENCE

By Louis H. Y. Chen[1] and Qi-Man Shao[2]

*National University of Singapore and University of Oregon*

We establish both uniform and nonuniform error bounds of the Berry–Esseen type in normal approximation under local dependence. These results are of an order close to the best possible if not best possible. They are more general or sharper than many existing ones in the literature. The proofs couple Stein's method with the concentration inequality approach.

**1. Introduction.** Since the work of Berry and Esseen in the 1940s, much has been done in the normal approximation for independent random variables. The standard tool has been the Fourier analytic method as developed by Esseen (1945). However, without independence, the Fourier analytic method becomes difficult to apply and bounds on the accuracy of approximation correspondingly difficult to find. In such situations, a method of Stein (1972) provides a much more viable alternative to the Fourier analytic method. Corresponding to calculating a Fourier transform and applying the inversion formula, it involves deriving a direct identity and solving a differential equation. As dependence is the rule rather than the exception in applications, Stein's method has become increasingly useful and important.

A crucial step in the Fourier analytic method for normal approximation is the use of a smoothing inequality originally due to Esseen (1945). The smoothing inequality is used to overcome the difficulty resulting from the nonsmoothness of the indicator function whose expectation is the distribution function. There is a correspondence of this in Stein's method, which is called the concentration inequality. It is originally due to Stein. Its simplest

Received May 2002; revised August 2003.
[1] Supported in part by Grant R-146-000-013-112 at National University of Singapore.
[2] Supported in part by NSF Grant DMS-01-03487, and Grants R-146-000-038-101 and R-1555-000-035-112 at National University of Singapore.

form is used in Ho and Chen (1978). More elaborate versions are proved in Chen (1986, 1998) and Chen and Shao (2001). In Chen and Shao (2001), it is developed also for obtaining nonuniform error bounds.

This paper is concerned with normal approximation under local dependence using Stein's method. Local dependence roughly means that certain subsets of the random variables are independent of those outside their respective "neighborhoods." No structure on the index set is assumed. Both uniform and nonuniform error bounds of the Berry–Esseen type are obtained and shown to be more general or sharper than many existing results in the literature. These include those of Shergin (1979), Prakasa Rao (1981), Baldi and Rinott (1989), Baldi, Rinott and Stein (1989), Rinott (1994) and Dembo and Rinott (1996).

The approach used in the paper is that of the concentration inequality. It is based on the ideas of Chen (1986) where the concentration inequality is derived differently from those in Chen and Shao (2001), due to the nonpositivity of the "covariance function." The uniform bounds obtained are improvements of those in Chen (1986), and the nonuniform bounds, which are proved by following the techniques in Chen and Shao (2001), are new in the literature. In proving the bounds, an attempt is made to achieve the best possible order for them. For example, the nonuniform bounds obtained are best possible as functions of the variables.

The forms of the bounds obtained are inspired by the results in Chen (1978), where necessary and sufficient conditions are proved for asymptotic normality of locally dependent random variables (termed finitely dependent in that paper). Such bounds deal successfully with those cases where the variance of a sum of $n$ random variables grows at a different rate from $n$. An example due to Erickson (1974) is used to illustrate this point.

This paper is organized as follows. The main results and their applications are given in Section 2. Two uniform and one nonuniform conditional concentration inequalities are proved in Section 3. The proofs of the uniform bounds are given in Section 4 and those of the nonuniform bounds in Section 5.

**2. Main results.** Throughout this paper let $\mathcal{J}$ be an index set and let $\{X_i, i \in \mathcal{J}\}$ be a random field with zero means and finite variances. Define $W = \sum_{i \in \mathcal{J}} X_i$ and assume that $\text{Var}(W) = 1$. Let $n$ be the cardinality of $\mathcal{J}$, let $F$ be the distribution function of $W$ and let $\Phi$ be the standard normal distribution function.

For $A \subset \mathcal{J}$, let $X_A$ denote $\{X_i, i \in A\}$, $A^c = \{j \in \mathcal{J} : j \notin A\}$, and let $|A|$ denote the cardinality of $A$. Adopt the notation: $a \wedge b = \min(a, b)$ and $a \vee b = \max(a, b)$.

We first introduce dependence assumptions and define notation that will be used throughout the paper.

(LD1) For each $i \in \mathcal{J}$, there exists $A_i \subset \mathcal{J}$ such that $X_i$ and $X_{A_i^c}$ are independent.
(LD2) For each $i \in \mathcal{J}$, there exist $A_i \subset B_i \subset \mathcal{J}$ such that $X_i$ is independent of $X_{A_i^c}$ and $X_{A_i}$ is independent of $X_{B_i^c}$.
(LD3) For each $i \in \mathcal{J}$, there exist $A_i \subset B_i \subset C_i \subset \mathcal{J}$ such that $X_i$ is independent of $X_{A_i^c}$, $X_{A_i}$ is independent of $X_{B_i^c}$ and $X_{B_i}$ is independent of $X_{C_i^c}$.
(LD4*) For each $i \in \mathcal{J}$, there exist $A_i \subset B_i \subset B_i^* \subset C_i^* \subset D_i^* \subset \mathcal{J}$ such that $X_i$ is independent of $X_{A_i^c}$, $X_{A_i}$ is independent of $X_{B_i^c}$, $X_{A_i}$ is independent of $\{X_{A_j}, j \in B_i^{*c}\}$, $\{X_{A_l}, l \in B_i^*\}$ is independent of $\{X_{A_j}, j \in C_i^{*c}\}$ and $\{X_{A_l}, l \in C_i^*\}$ is independent of $\{X_{A_j}, j \in D_i^{*c}\}$.

It is clear that (LD4*) implies (LD3), (LD3) yields (LD2) and (LD1) is the weakest assumption. Roughly speaking, (LD4*) is a version of (LD3) for $\{X_{A_i}, i \in \mathcal{J}\}$. On the other hand, in many cases (LD1) implies (LD2), (LD3) and (LD4*) with $B_i, C_i, B_i^*, C_i^*$ and $D_i^*$ defined as: $B_i = \bigcup_{j \in A_i} A_j$, $C_i = \bigcup_{j \in B_i} A_j$, $B_i^* = \bigcup_{j \in A_i} B_j$, $C_i^* = \bigcup_{j \in B_i^*} B_j$ and $D_i^* = \bigcup_{j \in C_i^*} B_j$. Other forms of local dependence have also been used in the literature, such as *dependency neighborhoods* in Rinott and Rotar (1996), where convergence rates of multivariate central limit theorem were obtained. Some of their results may not be covered by our theorems.

For each $i \in \mathcal{J}$, let $Y_i = \sum_{j \in A_i} X_j$. We define

$$\hat{K}_i(t) = X_i\{I(-Y_i \leq t < 0) - I(0 \leq t \leq -Y_i)\}, \quad K_i(t) = E\hat{K}_i(t),$$
(2.1) $$\hat{K}(t) = \sum_{i \in \mathcal{J}} \hat{K}_i(t), \quad K(t) = E\hat{K}(t) = \sum_{i \in \mathcal{J}} K_i(t).$$

Since $\mathrm{Var}(W) = 1$, we have

(2.2) $$\int_{-\infty}^{\infty} K(t)\,dt = EW^2 = 1.$$

2.1. *Uniform Berry–Esseen bounds.* The Berry–Esseen theorem [Berry (1941) and Esseen (1945); see, e.g., Petrov (1995)] states that if $\{X_i, i \in \mathcal{J}\}$ are independent with finite third moments, then there exists an absolute constant $C$ such that

$$\sup_{z \in R} |F(z) - \Phi(z)| \leq C \sum_{i \in \mathcal{J}} E|X_i|^3.$$

Here and throughout the paper, $C$ denotes an absolute constant which may have different values at different places. If, in addition, $X_i, i \in \mathcal{J}$, are identically distributed, then the bound is of the order $n^{-1/2}$, which is known to be the best possible.

The main objective of this section is to obtain general uniform Berry–Esseen bounds under various dependence assumptions with an aim to achieve the best possible orders. We first present a result under assumption (LD1).

THEOREM 2.1. *Under* (LD1), *we have*

$$\sup_z |F(z) - \Phi(z)| \leq r_1 + 4r_2 + 8r_3 + r_4 + 4.5r_5 + 1.5r_6, \tag{2.3}$$

*where*

$$r_1 = E\left|\sum_{i \in \mathcal{J}} (X_i Y_i - EX_i Y_i)\right|, \quad r_2 = \sum_{i \in \mathcal{J}} E|X_i Y_i| I(|Y_i| > 1),$$

$$r_3 = \sum_{i \in \mathcal{J}} E|X_i|(Y_i^2 \wedge 1), \quad r_4 = \sum_{i \in \mathcal{J}} E\{|WX_i|(Y_i^2 \wedge 1)\}, \tag{2.4}$$

$$r_5 = \int_{|t| \leq 1} \mathrm{Var}(\hat{K}(t)) \, dt, \quad r_6 = \left(\int_{|t| \leq 1} |t| \mathrm{Var}(\hat{K}(t)) \, dt\right)^{1/2}.$$

Since $r_1, r_2, r_3$ and $r_4$ depend on the moments of $\{X_i, Y_i, W\}$, they can be easily estimated (see Remark 2.1). The following alternative formulas of $r_5$ and $r_6$ may be useful. Let $\{X_i^*, i \in \mathcal{J}\}$ be an independent copy of $\{X_i, i \in \mathcal{J}\}$ and define $Y_i^* = \sum_{j \in A_i} X_j^*$. Then

$$\mathrm{Var}(\hat{K}(t)) = \sum_{i,j \in \mathcal{J}} E\{X_i X_j \{I(-Y_i \leq t < 0) - I(0 \leq t \leq -Y_i)\}$$

$$\times \{I(-Y_j \leq t < 0) - I(0 \leq t \leq -Y_j)\}$$

$$- X_i X_j^* \{I(-Y_i \leq t < 0) - I(0 \leq t \leq -Y_i)\}$$

$$\times \{I(-Y_j^* \leq t < 0) - I(0 \leq t \leq -Y_j^*)\}\},$$

and hence

$$r_5 = \sum_{i,j \in \mathcal{J}} E\{X_i X_j I(Y_i Y_j \geq 0)(|Y_i| \wedge |Y_j| \wedge 1)$$

$$- X_i X_j^* I(Y_i Y_j^* \geq 0)(|Y_i| \wedge |Y_j^*| \wedge 1)\}.$$

Similarly, we have

$$r_6^2 = \tfrac{1}{2} \sum_{i,j \in \mathcal{J}} E\{X_i X_j I(Y_i Y_j \geq 0)(|Y_i|^2 \wedge |Y_j|^2 \wedge 1)$$

$$- X_i X_j^* I(Y_i Y_j^* \geq 0)(|Y_i|^2 \wedge |Y_j^*|^2 \wedge 1)\}.$$

In particular, under assumption (LD2),

$$r_5 \leq \sum_{i,j \in \mathcal{J}, B_i B_j \neq \varnothing} E\{|X_i X_j|(|Y_i| \wedge |Y_j| \wedge 1) + |X_i X_j^*|(|Y_i| \wedge |Y_j^*| \wedge 1)\}$$

and

$$r_6^2 \leq \tfrac{1}{2} \sum_{i,j \in \mathcal{J}, B_i B_j \neq \varnothing} E\{|X_i X_j|(|Y_i|^2 \wedge |Y_j|^2 \wedge 1) + |X_i X_j^*|(|Y_i|^2 \wedge |Y_j^*|^2 \wedge 1)\}.$$

Thus, we have a much neater result under (LD2).

THEOREM 2.2. *Let $N(B_i) = \{j \in \mathcal{J} : B_j B_i \neq \varnothing\}$ and $2 < p \leq 4$. Assume that* (LD2) *is satisfied with $|N(B_i)| \leq \kappa$. Then*

$$\sup_z |F(z) - \Phi(z)| \leq (13 + 11\kappa) \sum_{i \in \mathcal{J}} (E|X_i|^{3 \wedge p} + E|Y_i|^{3 \wedge p})$$
$$(2.5) \qquad\qquad + 2.5 \left( \kappa \sum_{i \in \mathcal{J}} (E|X_i|^p + E|Y_i|^p) \right)^{1/2}.$$

*In particular, if $E|X_i|^p + E|Y_i|^p \leq \theta^p$ for some $\theta > 0$ and for each $i \in \mathcal{J}$, then*

$$(2.6) \qquad \sup_z |F(z) - \Phi(z)| \leq (13 + 11\kappa) n \theta^{3 \wedge p} + 2.5 \theta^{p/2} \sqrt{\kappa n},$$

*where $n = |\mathcal{J}|$.*

Note that in many cases $\kappa$ is bounded and $\theta$ is of order of $n^{-1/2}$. In those cases $\kappa n \theta^{3 \wedge p} + \theta^{p/2} \sqrt{\kappa n} = O(n^{-(p-2)/4})$, which is of the best possible order of $n^{-1/2}$ when $p = 4$. However, the cost is the existence of fourth moments. To reduce the assumption on moments, we need a stronger condition.

THEOREM 2.3. *Suppose that* (LD3) *is satisfied. Let $N(C_i) = \{j \in \mathcal{J} : C_i B_j \neq \varnothing\}$,*

$$\tilde{W}_i = \sum_{j \in N(C_i)^c} X_j, \qquad \sigma_i^2 = \operatorname{Var}(\tilde{W}_i) \quad \text{and} \quad \lambda = 1 \vee \max_{i \in \mathcal{J}} (1/\sigma_i).$$

*Then*

$$(2.7) \quad \sup_z |F(z) - \Phi(z)| \leq 4\lambda^{3/2} (r_2 + r_3 + r_7 + r_8 + r_9 + r_{10} + r_{11} + r_{12}),$$

*where $Z_i = \sum_{j \in B_i} X_j$,*

$$r_7 = \sum_{i \in \mathcal{J}} E\{|X_i Y_i| I(|X_i| > 1)\},$$
$$r_8 = \sum_{i \in \mathcal{J}} E\{|X_i| I(|X_i| \leq 1)(|Y_i| \wedge 1)|Z_i|\},$$
$$r_9 = \sum_{i \in \mathcal{J}} E\{|W X_i| I(|X_i| \leq 1)(|Y_i| \wedge 1)(|Z_i| \wedge 1)\},$$
$$(2.8) \quad r_{10} = \sum_{i,j \in \mathcal{J}, B_i B_j \neq \varnothing} E\{|X_i X_j|(|Y_i| \wedge |Y_j| \wedge 1) + |X_i X_j^*|(|Y_i| \wedge |Y_j^*| \wedge 1)\},$$
$$r_{11} = \sum_{i \in \mathcal{J}} P(|X_i| > 1) E|X_i|(|Y_i| \wedge 1),$$
$$r_{12} = \sum_{i \in \mathcal{J}} E\{(|W| + 1)(|Z_i| \wedge 1)\} E|X_i|(|Y_i| \wedge 1),$$

*and $(X_i^*, Y_i^*)$ is an independent copy of $(X_i, Y_i)$.*

In particular, we have:

THEOREM 2.4. *Let $2 < p \leq 3$. Assume that* (LD3) *is satisfied with* $\max(|N(C_i)|, |\{j : i \in C_j\}|) \leq \kappa$. *Then*

$$\sup_z |F(z) - \Phi(z)| \leq 75\kappa^{p-1} \sum_{i \in \mathcal{J}} E|X_i|^p. \tag{2.9}$$

Rinott (1994) and Dembo and Rinott (1996) obtained uniform bounds of order $n^{-1/2}$ when $X_i$ is bounded with order of $n^{-1/2}$ under a different local dependence assumption which appears to be weaker than (LD2). However, their approach does not seem to be extendable to random variables which are not necessarily bounded.

REMARK 2.1. Although $r_4$ involves $W$, there are several ways to bound it. When there is no additional assumption besides (LD1), we can use the following estimate:

$$E|WX_i|(Y_i^2 \wedge 1)$$

$$\leq E|WX_i|\left(|A_i| \sum_{j \in A_i} (X_j^2 \wedge 1)\right)$$

$$\leq |A_i| \sum_{j \in A_i} E|WX_i|(X_j^2 \wedge 1) I(|X_i| \leq |X_j|)$$

$$\quad + |A_i| \sum_{j \in A_i} E|WX_i|(X_j^2 \wedge 1) I(|X_i| > |X_j|)$$

$$\leq |A_i| \sum_{j \in A_i} E|W||X_j|(X_j^2 \wedge 1) + |A_i| \sum_{j \in A_i} E|W||X_i|(X_i^2 \wedge 1)$$

$$\leq |A_i| \sum_{j \in A_i} E|W - Y_j| E|X_j|(X_j^2 \wedge 1) + |A_i| \sum_{j \in A_i} E|Y_j X_j|(|X_j| \wedge 1)$$

$$\quad + |A_i|^2 E|W - Y_i| E|X_i|(X_i^2 \wedge 1) + |A_i|^2 E|Y_i X_i|(|X_i| \wedge 1)$$

$$\leq |A_i| \sum_{j \in A_i} (1 + E|Y_j|) E|X_j|(X_j^2 \wedge 1) + |A_i| \sum_{j \in A_i} E|Y_j X_j|(|X_j| \wedge 1)$$

$$\quad + |A_i|^2 (1 + E|Y_i|) E|X_i|(X_i^2 \wedge 1) + |A_i|^2 E|Y_i X_i|(|X_i| \wedge 1).$$

2.2. *Nonuniform Berry–Esseen bound.* Nonuniform bounds were first obtained by Esseen (1945) for independent and identically distributed random variables $\{X_i, i \in \mathcal{J}\}$. These were improved to $CnE|X_1|^3/(1+|x|^3)$ by Nagaev (1965). Bikelis (1966) generalized Nagaev's result to

$$|F(z) - \Phi(z)| \leq \frac{C \sum_{i \in \mathcal{J}} E|X_i|^3}{1 + |z|^3}$$

for independent and not necessarily identically distributed random variables.

In this section we present a general nonuniform bound for locally dependent random fields $\{X_i, i \in \mathcal{J}\}$ under (LD4*).

THEOREM 2.5. *Assume that $E|X_i|^p < \infty$ for $2 < p \leq 3$ and that* (LD4*) *is satisfied. Let $\kappa = \max_{i \in \mathcal{J}} \max(|D_i^*|, |\{j : i \in D_j^*\}|)$. Then*

$$(2.10) \qquad |F(z) - \Phi(z)| \leq C\kappa^p(1+|z|)^{-p} \sum_{i \in \mathcal{J}} E|X_i|^p.$$

2.3. *m-dependent random fields.* Let $d \geq 1$ and let $Z^d$ denote the $d$-dimensional space of positive integers. The distance between two points $i = (i_1, \ldots, i_d)$ and $j = (j_1, \ldots, j_d)$ in $Z^d$ is defined by $|i-j| = \max_{1 \leq l \leq d} |i_l - j_l|$ and the distance between two subsets $A$ and $B$ of $Z^d$ is defined by $\rho(A,B) = \inf\{|i-j| : i \in A, j \in B\}$. For a given subset $\mathcal{J}$ of $Z^d$, a set of random variables $\{X_i, i \in \mathcal{J}\}$ is said to be an $m$-dependent random field if $\{X_i, i \in A\}$ and $\{X_j, j \in B\}$ are independent whenever $\rho(A,B) > m$, for any subsets $A$ and $B$ of $\mathcal{J}$.

Thus choosing $A_i = \{j : |j-i| \leq m\} \cap \mathcal{J}, B_i = \{j : |j-i| \leq 2m\} \cap \mathcal{J}, C_i = \{j : |j-i| \leq 3m\} \cap \mathcal{J}, B_i^* = \{j : |j-i| \leq 3m\} \cap \mathcal{J}, C_i^* = \{j : |j-i| \leq 6m\} \cap \mathcal{J}$ and $D_i^* = \{j : |j-i| \leq 9m\} \cap \mathcal{J}$ in Theorems 2.4 and 2.5 yields a uniform and a nonuniform bound.

THEOREM 2.6. *Let $\{X_i, i \in \mathcal{J}\}$ be an $m$-dependent random field with zero means and finite $E|X_i|^p < \infty$ for $2 < p \leq 3$. Then*

$$(2.11) \qquad \sup_z |F(z) - \Phi(z)| \leq 75(10m+1)^{(p-1)d} \sum_{i \in \mathcal{J}} E|X_i|^p$$

*and*

$$(2.12) \quad |F(z) - \Phi(z)| \leq C(1+|z|)^{-p}(19)^{pd}(m+1)^{(p-1)d} \sum_{i \in \mathcal{J}} E|X_i|^p.$$

Here we have reduced the $m$-dependent random field to a one-dependent random field by taking blocks and then applied (2.10) to get (2.12). The result (2.11) was previously obtained by Shergin (1979) without specifying the absolute constant. For nonuniform bounds, results weaker than (2.12) have been obtained in the literature. See, for example, Prakasa Rao (1981) and Heinrich (1984). However, the result in Prakasa Rao (1981) is far from best possible even for independent random fields, while Heinrich (1984) is the best possible only for the i.i.d. case. For other uniform and nonuniform Berry–Esseen bounds for $m$-dependent and weakly dependent random variables, see Tihomirov (1980), Dasgupta (1992) and Sunklodas (1999). In Sunklodas (1999) a lower bound is also given.

2.4. *Examples.* In this section we give three examples discussed in literature to illustrate the usefulness of our general results.

2.4.1. *Graph dependency.* This example was discussed in Baldi and Rinott (1989) and Rinott (1994), where some results on uniform bound were obtained.

Consider a set of random variables $\{X_i, i \in \mathcal{V}\}$ indexed by the vertices of a graph $\mathcal{G} = (\mathcal{V}, \mathcal{E})$. $\mathcal{G}$ is said to be a dependency graph if, for any pair of disjoint sets $\Gamma_1$ and $\Gamma_2$ in $\mathcal{V}$ such that no edge in $\mathcal{E}$ has one endpoint in $\Gamma_1$ and the other in $\Gamma_2$, the sets of random variables $\{X_i, i \in \Gamma_1\}$ and $\{X_i, i \in \Gamma_2\}$ are independent. Let $D$ denote the maximal degree of $G$, that is, the maximal number of edges incident to a single vertex. Let $A_i = \{j \in \mathcal{V}: \text{there is an edge connecting } j \text{ and } i\}$, $B_i = \bigcup_{j \in A_i} A_j$, $C_i = \bigcup_{j \in B_i} A_j$, $B_i^* = \bigcup_{j \in A_i} B_j$, $C_i^* = \bigcup_{j \in B_i^*} B_j$ and $D_i^* = \bigcup_{j \in C_i^*} B_j$. An application of Theorems 2.4 and 2.5 yields the following theorem.

THEOREM 2.7. *Let $\{X_i, i \in \mathcal{V}\}$ be random variables indexed by the vertices of a dependency graph. Put $W = \sum_{i \in \mathcal{V}} X_i$. Assume that $EW^2 = 1$, $EX_i = 0$ and $E|X_i|^p \leq \theta^p$ for $i \in \mathcal{V}$ and for some $\theta > 0$. Then*

$$\sup_z |P(W \leq z) - \Phi(z)| \leq 75 D^{5(p-1)} |\mathcal{V}| \theta^p \tag{2.13}$$

*and for $z \in R$,*

$$|P(W \leq z) - \Phi(z)| \leq C(1 + |z|)^{-p} D^{5p} |\mathcal{V}| \theta^p. \tag{2.14}$$

While (2.13) compares favorably with those of Baldi and Rinott (1989), (2.14) is new.

2.4.2. *The number of local maxima on a graph.* Consider a graph $\mathcal{G} = (\mathcal{V}, \mathcal{E})$ (which is not necessarily a dependency graph) and independently identically distributed continuous random variables $\{Y_i, i \in \mathcal{V}\}$. For $i \in \mathcal{V}$, define the 0–1 indicator variable

$$X_i = \begin{cases} 1, & \text{if } Y_i > Y_j \text{ for all } j \in N_i, \\ 0, & \text{otherwise,} \end{cases}$$

where $N_i = \{j \in \mathcal{V}: d(i,j) = 1\}$ and $d(i,j)$ denotes the shortest path distance between the vertices $i$ and $j$. Note that $d(i,j) = 1$ iff $i$ and $j$ are neighbors, so $X_i = 1$ indicates that $Y_i$ is a local maximum. Let $W = \sum_{i \in \mathcal{V}} X_i$ be the number of local maxima. If $(\mathcal{V}, \mathcal{E})$ is a regular graph, that is, all vertices have the same degree $d$, then by Baldi, Rinott and Stein (1989), $EW = |\mathcal{V}|/(d+1)$,

$$\sigma^2 = \text{Var}(W) = \sum_{i,j \in \mathcal{V}, d(i,j)=2} s(i,j)(2d+2-s(i,j))^{-1}(d+1)^{-2}$$

and
$$\sup_z |P(W \leq z) - \Phi((z - EW)/\sigma)| \leq C\sigma^{-1/2},$$

where $s(i,j) = |N_i \cap N_j|$.

Theorem 2.8 is obtained by applying Theorem 2.2. The uniform bound, which improves $\sigma^{-1/2}$ of Baldi, Rinott and Stein (1989) to $\sigma^{-1}$, is similar to that of Dembo and Rinott (1996). However, the nonuniform bound is new.

THEOREM 2.8. *We have*

(2.15) $$\sup_z |P(W \leq z) - \Phi((z - EW)/\sigma)| \leq Cd^2|\mathcal{V}|/\sigma^3$$

*and*

(2.16) $$|P(W \leq z) - \Phi((z - EW)/\sigma)| \leq C(1 + |z|)^{-3}d^5|\mathcal{V}|/\sigma^3.$$

THEOREM 2.9. *We have*

(2.17) $$|P(W \leq z) - \Phi((z - np)/\sigma)| \leq \frac{Cr^2}{(1 + |z|^3)\sqrt{np(1-p)}}.$$

2.4.3. *One-dependence with $o(n)$ variance.* This example was discussed in Erickson (1974). Define a sequence of bounded, symmetric and identically distributed random variables $X_1, \ldots, X_n$ with $EX_i^2 = 1$ as follows. Let $X_1$ and $X_2$ be independent bounded and symmetric random variables with variance 1 and put $B_k^2 = \text{Var}(\sum_{i=1}^k X_i)$. For $k \geq 2$, define $X_{k+1} = -X_k$ if $B_k^2 > k^{1/2}$ and define $X_{k+1}$ to be independent of $X_1, \ldots, X_k$ if $B_k^2 \leq k^{1/2}$. It is clear that $X_1, \ldots, X_n$ is a one-dependent sequence, $|B_n^2 - n^{1/2}| \leq 2$ and $\sum_{i=1}^n X_i$ is a sum of $B_n^2 \sim n^{1/2}$ terms of independent and identically distributed random variables. By the Berry–Esseen theorem,
$$\sup_z |P(W \leq z) - \Phi(z/B_n)| \leq Cn^{-1/4}E|X_1|^3$$

and the order $n^{-1/4}$ is correct. While the bound in (2.11) generalizes and improves many others, it is asymptotically $Cn^{1/4}E|X_1|^3$, which goes to $\infty$. On the other hand, Theorem 2.3 gives the correct order $n^{-1/4}$. To see this, we observe that if $X_i$ is independent of all other random variables, then we can choose $A_i = B_i = C_i = \{i\}$ and $X_i = Y_i = Z_i$. On the other hand, if $X_i = -X_{i-1}$ or $-X_{i+1}$, then $A_i = B_i = C_i = \{i - 1, i\}$ or $\{i, i + 1\}$, respectively. In this case $Y_i = Z_i = 0$ identically. Consequently, the right-hand side of (2.7) is bounded by

$$CE|X_1/B_n|^3 \times \text{number of } X_i \text{ which is independent}$$
$$\text{of all the other random variables}$$
$$\leq CE|X_1/B_n|^3 B_n^2 = CE|X_1|^3 B_n^{-1} \leq Cn^{-1/4}E|X_1|^3,$$

as desired.

**3. Concentration inequalities.** The concentration inequality in normal approximation using Stein's method plays a role corresponding to that of the smoothing inequality of Esseen (1945) in the Fourier analytic method. It is used to overcome the difficulty caused by the nonsmoothness of the indicator function whose expectation is the distribution function. In this section we establish two uniform and one nonuniform conditional concentration inequalities. We first prove Propositions 3.1, 3.2 and 3.3, which will be used in the proofs of Theorems 2.1, 2.2 and 2.5, respectively. Esseen (1968), Petrov (1995) and others have obtained many uniform concentration inequalities for sums of independent random variables. Our uniform concentration inequalities are different from theirs except in the i.i.d. case.

Let $\{X_i, i \in \mathcal{J}\}$ be a random field with $EX_i = 0$ and $EX_i^2 < \infty$. Put $W = \sum_{i \in \mathcal{J}} X_i$. Assume that $EW^2 = 1$.

PROPOSITION 3.1. *Assume* (LD1). *Then for any real numbers $a < b$,*

$$(3.1) \qquad P(a \leq W \leq b) \leq 0.625(b-a) + 4r_2 + 2.125r_3 + 4r_5,$$

*where $r_2$, $r_3$ and $r_5$ are as defined in* (2.4).

PROOF. Let $\alpha = r_3$ and define

$$(3.2) \quad f(w) = \begin{cases} -\dfrac{b-a+\alpha}{2}, & \text{for } w \leq a - \alpha, \\ \dfrac{1}{2\alpha}(w-a+\alpha)^2 - \dfrac{b-a+\alpha}{2}, & \text{for } a-\alpha < w \leq a, \\ w - \dfrac{a+b}{2}, & \text{for } a < w \leq b, \\ -\dfrac{1}{2\alpha}(w-b-\alpha)^2 + \dfrac{b-a+\alpha}{2}, & \text{for } b < w \leq b+\alpha, \\ \dfrac{b-a+\alpha}{2}, & \text{for } w > b+\alpha. \end{cases}$$

Then $f'$ is a continuous function given by

$$f'(w) = \begin{cases} 1, & \text{for } a \leq w \leq b, \\ 0, & \text{for } w \leq a-\alpha \text{ or } w \geq b+\alpha, \\ \text{linear}, & \text{for } a-\alpha \leq w \leq a \text{ or } b \leq w \leq b+\alpha. \end{cases}$$

Clearly, $|f(w)| \leq (b-a+\alpha)/2$. With this $f$, $Y_i$, and $\hat{K}(t)$ and $K(t)$ as defined in (2.1), we have

$$(b-a+\alpha)/2 \geq EWf(W) = \sum_{i \in \mathcal{J}} E\{X_i(f(W) - f(W-Y_i))\}$$

$$= \sum_{i \in \mathcal{J}} E\left\{X_i \int_{-Y_i}^0 f'(W+t)\,dt\right\}$$

(3.3)
$$= \sum_{i \in \mathcal{J}} E\left\{\int_{-\infty}^{\infty} f'(W+t)\hat{K}_i(t)\,dt\right\}$$
$$= E\int_{-\infty}^{\infty} f'(W+t)\hat{K}(t)\,dt := H_1 + H_2 + H_3 + H_4,$$

where
$$H_1 = Ef'(W)\int_{|t|\leq 1} K(t)\,dt,$$
$$H_2 = E\left\{\int_{|t|\leq 1}(f'(W+t) - f'(W))K(t)\,dt\right\},$$
$$H_3 = E\left\{\int_{|t|>1} f'(W+t)\hat{K}(t)\,dt\right\},$$
$$H_4 = E\left\{\int_{|t|\leq 1} f'(W+t)(\hat{K}(t) - K(t))\,dt\right\}.$$

Clearly, by (2.2),

(3.4)
$$H_1 = Ef'(W)\left\{1 - \int_{|t|>1} K(t)\,dt\right\}$$
$$\geq Ef'(W)(1 - r_2) \geq P(a \leq W \leq b) - r_2$$

and

(3.5)
$$|H_3| \leq \sum_{i \in \mathcal{J}} E|X_iY_i|I(|Y_i| > 1) = r_2.$$

By the Cauchy inequality,

(3.6)
$$|H_4| \leq \frac{1}{8}E\int_{|t|\leq 1}[f'(W+t)]^2\,dt + 2E\int_{|t|\leq 1}(\hat{K}(t) - K(t))^2\,dt$$
$$\leq \frac{b - a + 2\alpha}{8} + 2r_5.$$

To bound $H_2$, let
$$L(\alpha) = \sup_{x \in R} P(x \leq W \leq x + \alpha).$$

Then by writing
$$H_2 = E\int_0^1\int_0^t f''(W+s)\,ds\,K(t)\,dt - E\int_{-1}^0\int_t^0 f''(W+s)\,ds\,K(t)\,dt$$
$$= \alpha^{-1}\int_0^1\int_0^t \{P(a-\alpha \leq W+s \leq a) - P(b \leq W+s \leq b+\alpha)\}\,ds\,K(t)\,dt$$
$$- \alpha^{-1}\int_{-1}^0\int_t^0 \{P(a-\alpha \leq W+s \leq a)$$
$$- P(b \leq W+s \leq b+\alpha)\}\,ds\,K(t)\,dt,$$

we have

$$
\begin{aligned}
|H_2| &\leq \alpha^{-1} \int_0^1 \int_0^t L(\alpha) \, ds \, |K(t)| \, dt + \alpha^{-1} \int_{-1}^0 \int_t^0 L(\alpha) \, ds \, |K(t)| \, dt \\
&= \alpha^{-1} L(\alpha) \int_{|t| \leq 1} |tK(t)| \, dt \leq 0.5 \alpha^{-1} L(\alpha) r_3 = 0.5 L(\alpha).
\end{aligned}
\tag{3.7}
$$

It follows from (3.3)–(3.7) that

$$
P(a \leq W \leq b) \leq 0.625(b-a) + 0.75\alpha + 2r_2 + 2r_5 + 0.5L(\alpha). \tag{3.8}
$$

Substituting $a = x$ and $b = x + \alpha$ in (3.8), we obtain

$$
L(\alpha) \leq 1.375\alpha + 2r_2 + 2r_5 + 0.5L(\alpha)
$$

and hence

$$
L(\alpha) \leq 2.75\alpha + 4r_2 + 4r_5. \tag{3.9}
$$

Finally, combining (3.8) and (3.9), we obtain (3.1). □

PROPOSITION 3.2. *Assume* (LD3). *Let* $\xi = \xi_i = (X_i, Y_i, Z_i)$, *where* $Y_i = \sum_{j \in A_i} X_j$ *and* $Z_i = \sum_{j \in B_i} X_j$. *For Borel measurable functions* $a_\xi$ *and* $b_\xi$ *of* $\xi$ *such that* $a_\xi \leq b_\xi$, *we have*

$$
\begin{aligned}
&P^\xi(a_\xi \leq W \leq b_\xi) \\
&\quad \leq 0.625 \sigma_i^{-1}(b_\xi - a_\xi) + 4\sigma_i^{-2} r_2 + 2.125 \sigma_i^{-3} r_3 + 4\sigma_i^{-3} r_{10} \quad \text{a.s.},
\end{aligned}
\tag{3.10}
$$

*where* $\sigma_i$ *and* $r_{10}$ *are as defined in Theorem* 2.3, *and* $P^\xi(\cdot)$ *denotes the conditional probability given* $\xi$.

PROOF. We use the same notation as in Theorem 2.3 and follow the same line of the proof as that of Proposition 3.1. Let $f_\xi$ be defined similarly as in (3.2) such that $f_\xi((a_\xi + b_\xi)/2) = 0$ and $f'_\xi$ is a continuous function given by

$$
f'_\xi(w) = \begin{cases} 1, & \text{for } a_\xi \leq w \leq b_\xi, \\ 0, & \text{for } w \leq a_\xi - \alpha \text{ or } w \geq b_\xi + \alpha, \\ \text{linear}, & \text{for } a_\xi - \alpha \leq w \leq a_\xi \text{ or } b_\xi \leq w \leq b_\xi + \alpha, \end{cases}
$$

where $\alpha = \sigma_i^{-2} r_3$. Then, $|f_\xi(w)| \leq (b_\xi - a_\xi + \alpha)/2$. Put

$$
\hat{M}(t) = \sum_{j \in N(C_i)^c} \hat{K}_j(t) \quad \text{and} \quad M(t) = \sum_{j \in N(C_i)^c} K_j(t).
$$

Observe that $X_j$ and $\{\xi, W - Y_j\}$ are independent for $j \in N(C_i)^c$. Similarly to (3.3),

$$
\begin{aligned}
0.5(b_\xi - a_\xi + \alpha)\sigma_i &\geq E\{\tilde{W}_i f_\xi(W)\} \\
&= \sum_{j \in N(C_i)^c} E^\xi\{X_j(f_\xi(W) - f_\xi(W - Y_j))\} \\
&= E^\xi\left\{\int_{-\infty}^{\infty} f'_\xi(W + t)\hat{M}(t)\,dt\right\} \\
&:= H_{1,\xi} + H_{2,\xi} + H_{3,\xi} + H_{4,\xi},
\end{aligned}
\tag{3.11}
$$

where $E^\xi(\cdot)$ denotes the conditional expectation given $\xi$,

$$H_{1,\xi} = E^\xi f'_\xi(W) \int_{|t|\leq 1} M(t)\,dt,$$

$$H_{2,\xi} = E^\xi\left\{\int_{|t|\leq 1}(f'_\xi(W+t) - f'_\xi(W))M(t)\,dt\right\},$$

$$H_{3,\xi} = E^\xi\left\{\int_{|t|>1} f'_\xi(W+t)\hat{M}(t)\,dt\right\},$$

$$H_{4,\xi} = E^\xi\left\{\int_{|t|\leq 1} f'_\xi(W+t)(\hat{M}(t) - M(t))\,dt\right\}.$$

Note that $\xi$ and $\hat{M}(t)$ are independent. Analogously to (3.4)–(3.6),

$$H_{1,\xi} \geq P^\xi(a_\xi \leq W \leq b_\xi)\sigma_i^2 - r_2, \tag{3.12}$$

$$
\begin{aligned}
|H_{3,\xi}| &\leq \sum_{j \in N(C_i)^c} E^\xi\{|X_j Y_j| I(|Y_j| > 1)\} \\
&= \sum_{j \in N(C_i)^c} E\{|X_j Y_j| I(|Y_j| > 1)\} \leq r_2,
\end{aligned}
\tag{3.13}
$$

$$
\begin{aligned}
|H_{4,\xi}| &\leq 0.125\sigma_i(b_\xi - a_\xi + 2\alpha) + 2\sigma_i^{-1}E\int_{|t|\leq 1}(\hat{M}(t) - M(t))^2\,dt \\
&= 0.125\sigma_i(b_\xi - a_\xi + 2\alpha) + 2\sigma_i^{-1}\rho,
\end{aligned}
\tag{3.14}
$$

where $\rho = \int_{|t|\leq 1} \text{Var}(\hat{M}(t))\,dt$.

To bound $H_{2,\xi}$, define

$$L_\xi(\alpha) = \lim_{k\to\infty}\sup_{x\in Q} P^\xi(x - 1/k \leq W \leq x + 1/k + \alpha),$$

where $Q$ is the set of rational numbers and, with a little abuse of notation, we regard $P^\xi$ as a regular conditional probability given $\xi$. Following the proof of (3.7) yields

$$(3.15)\quad |H_{2,\xi}| \leq \alpha^{-1}L_\xi(\alpha)\int_{|t|\leq 1} |tM(t)|\,dt \leq 0.5\alpha^{-1}L_\xi(\alpha)r_3 = 0.5\sigma_i^2 L_\xi(\alpha).$$

Thus by (3.12)–(3.15),

$$
\begin{aligned}
(3.16) \quad & P^\xi(a_\xi \leq W \leq b_\xi) \\
& \leq 0.625\sigma_i^{-1}(b_\xi - a_\xi) + 0.75\alpha\sigma_i^{-1} + 2\sigma_i^{-2}r_2 + 2\sigma_i^{-3}\rho + 0.5L_\xi(\alpha).
\end{aligned}
$$

Now substitute $a_\xi = x - 1/k$ and $b_\xi = x + 1/k + \alpha$ in (3.16). By taking the supremum over $x \in Q$ and then letting $k \to \infty$, we obtain

$$L_\xi(\alpha) \leq 1.375\alpha\sigma_i^{-1} + 2\sigma_i^{-2}r_2 + 2\sigma_i^{-3}\rho + 0.5L_\xi(\alpha)$$

and hence

$$(3.17) \qquad L_\xi(\alpha) \leq 2.75\alpha\sigma_i^{-1} + 4\sigma_i^{-2}r_2 + 4\sigma_i^{-3}\rho.$$

Combining (3.16) and (3.17), we obtain

$$
\begin{aligned}
(3.18) \quad & P^\xi(a_\xi \leq W \leq b_\xi) \\
& \leq 0.625\sigma_i^{-1}(b_\xi - a_\xi) + 4\sigma_i^{-2}r_2 + 2.125\sigma_i^{-3}r_3 + 4\sigma_i^{-3}\rho.
\end{aligned}
$$

It remains to prove that $\rho \leq r_{10}$. Let $(X_j^*, Y_j^*)$ be an independent copy of $(X_j, Y_j)$. Note that $\hat{K}_j(t)$ and $\hat{K}_l(t)$ are independent for $l \in N(B_j)^c$. Direct computations yield

$$
\begin{aligned}
\rho &= \sum_{j \in N(C_i)^c} \sum_{l \in N(C_i)^c} \int_{|t| \leq 1} \operatorname{Cov}(\hat{K}_j(t), \hat{K}_l(t))\,dt \\
&= \sum_{j \in N(C_i)^c} \sum_{l \in N(B_j)N(C_i)^c} \int_{|t| \leq 1} \operatorname{Cov}(\hat{K}_j(t), \hat{K}_l(t))\,dt \\
&= \sum_{j \in N(C_i)^c} \sum_{l \in N(B_j)N(C_i)^c} E\{X_l X_j I(Y_l Y_j \geq 0)(|Y_l| \wedge |Y_j| \wedge 1) \\
&\qquad\qquad - X_l X_j^* I(Y_l Y_j^* \geq 0)(|Y_l| \wedge |Y_j^*| \wedge 1)\} \\
&\leq \sum_{j \in \mathcal{J}} \sum_{l \in N(B_j)} E\{|X_l X_j|(|Y_l| \wedge |Y_j| \wedge 1) + |X_l X_j^*|(|Y_l| \wedge |Y_j^*| \wedge 1)\} \\
&= r_{10}.
\end{aligned}
$$

This completes the proof of the proposition. $\square$

For obtaining a nonuniform conditional concentration inequality, we need two lemmas on moment inequalities for locally dependent random fields.

LEMMA 3.1. *Let $\{X_i, i \in \mathcal{J}\}$ be a random field satisfying* (LD3) *and let $\xi_i$ be a measurable function of $X_i$ with $E\xi_i = 0$ and $E\xi_i^4 < \infty$ for each $i \in \mathcal{J}$. Let $T = \sum_{i \in \mathcal{J}} \xi_i$ and $\sigma^2 = E(T^2)$. Then*

$$(3.19) \qquad \sigma^2 = \sum_{i \in \mathcal{J}} \sum_{j \in A_i} E\xi_i \xi_j$$

and, for $a > 0$,

$$(3.20) \quad ET^4 = 3\sigma^4 - 6\sum_{i \in \mathcal{J}} E(\xi_i \xi_{A_i}) E(\xi_{B_i} \xi_{C_i}) + 3\sum_{i \in \mathcal{J}} E(\xi_i \xi_{A_i}) E\xi_{B_i}^2$$

$$- 3\sum_{i \in \mathcal{J}} E(\xi_i \xi_{A_i}^2 \xi_{B_i}) + \sum_{i \in \mathcal{J}} E(\xi_i \xi_{A_i}^3)$$

$$+ 6\sum_{i \in \mathcal{J}} E(\xi_i \xi_{A_i} \xi_{B_i} \xi_{C_i}) - 3\sum_{i \in \mathcal{J}} E(\xi_i \xi_{A_i} \xi_{B_i}^2)$$

$$(3.21) \quad \leq 3\sigma^4 + 5.5 \sum_{i \in \mathcal{J}} \{a^3 E\xi_i^4 + a^{-1} E\xi_{A_i}^4 + a^{-1} E\xi_{B_i}^4 + a^{-1} E\xi_{C_i}^4\},$$

where

$$\xi_{A_i} = \sum_{j \in A_i} \xi_j, \qquad \xi_{B_i} = \sum_{j \in B_i} \xi_j, \qquad \xi_{C_i} = \sum_{j \in C_i} \xi_j.$$

In particular, we have

$$(3.22) \quad \sigma^2 \leq \kappa_1 \sum_{i \in \mathcal{J}} E\xi_i^2$$

and

$$(3.23) \quad ET^4 \leq 3\sigma^4 + 22\kappa_1^3 \sum_{i \in \mathcal{J}} E\xi_i^4,$$

where $\kappa_1 = \max_{i \in J} \max(|C_i|, |\{j \in \mathcal{J} : i \in C_j\}|)$.

PROOF. By (LD3), $\xi_i$ and $\{\xi_j, j \in A_i^c\}$ are independent and this implies (3.19). Note that for each $i \in \mathcal{J}$, $\xi_i$ and $T - \xi_{A_i}$ are independent, $\{\xi_i, \xi_{A_i}\}$ and $T - \xi_{B_i}$ are independent and $\{\xi_i, \xi_{A_i}, \xi_{B_i}\}$ and $T - \xi_{C_i}$ are independent. Therefore

$$ET^4 = \sum_{i \in \mathcal{J}} E\{\xi_i (T^3 - (T - \xi_{A_i})^3)\}$$

$$= \sum_{i \in \mathcal{J}} 3E\{\xi_i \xi_{A_i} T^2\} - 3E\{\xi_i \xi_{A_i}^2 T\} + E\{\xi_i \xi_{A_i}^3\}$$

$$(3.24) \quad = 3\sum_{i \in \mathcal{J}} E\{\xi_i \xi_{A_i}\} E(T - \xi_{B_i})^2 + 3\sum_{i \in \mathcal{J}} E\{\xi_i \xi_{A_i}(T^2 - (T - \xi_{B_i})^2)\}$$

$$- 3\sum_{i \in \mathcal{J}} E\{\xi_i \xi_{A_i}^2 \xi_{B_i}\} + \sum_{i \in \mathcal{J}} E\{\xi_i \xi_{A_i}^3\}$$

$$= 3\sum_{i \in \mathcal{J}} E\{\xi_i \xi_{A_i}\}(\sigma^2 - 2E\xi_{B_i}\xi_{C_i} + E\xi_{B_i}^2)$$

$$+ 3\sum_{i \in \mathcal{J}} E\{\xi_i \xi_{A_i}(2\xi_{B_i}\xi_{C_i} - \xi_{B_i}^2)\}$$

$$- 3 \sum_{i \in \mathcal{J}} E\{\xi_i \xi_{A_i}^2 \xi_{B_i}\} + \sum_{i \in \mathcal{J}} E\{\xi_i \xi_{A_i}^3\}$$

$$= 3\sigma^4 - 6 \sum_{i \in \mathcal{J}} E(\xi_i \xi_{A_i}) E(\xi_{B_i} \xi_{C_i}) + 3 \sum_{i \in \mathcal{J}} E(\xi_i \xi_{A_i}) E\xi_{B_i}^2$$

$$- 3 \sum_{i \in \mathcal{J}} E(\xi_i \xi_{A_i}^2 \xi_{B_i}) + \sum_{i \in \mathcal{J}} E(\xi_i \xi_{A_i}^3)$$

$$+ 6 \sum_{i \in \mathcal{J}} E(\xi_i \xi_{A_i} \xi_{B_i} \xi_{C_i}) - 3 \sum_{i \in \mathcal{J}} E(\xi_i \xi_{A_i} \xi_{B_i}^2).$$

Now the Cauchy inequality implies that, for any $a > 0$ and any random variables $u_1, u_2, u_3, u_4$,

(3.25) $\quad E|u_1 u_2 u_3 u_4|$
$\leq \frac{1}{4}\{a^3 E|u_1|^4 + a^{-1} E|u_2|^4 + a^{-1} E|u_3|^4 + a^{-1} E|u_4|^4\}.$

It follows that the right-hand side of (3.24) is bounded by $5.5 \sum_{i \in \mathcal{J}} \{a^3 E\xi_i^4 + a^{-1} E\xi_{A_i}^4 + a^{-1} E\xi_{B_i}^4 + a^{-1} E\xi_{C_i}^4\}$. This proves (3.21).

To prove (3.22) and (3.23), put $A_i^{-1} = \{j \in \mathcal{J} : i \in A_j\}$ and $C_i^{-1} = \{j \in \mathcal{J} : i \in C_j\}$. Then $\kappa_1 = \max_{i \in \mathcal{J}}(|C_i| \vee |C_i^{-1}|)$. By the $C_r$ inequality and (3.19),

$$\sigma^2 \leq \sum_{i \in \mathcal{J}} \sum_{j \in A_i} 0.5\{E\xi_i^2 + E\xi_j^2\}$$

$$\leq 0.5 \kappa_1 \sum_{i \in \mathcal{J}} E\xi_i^2 + 0.5 \sum_{j \in \mathcal{J}} \sum_{i \in C_j^{-1}} E\xi_j^2 \leq \kappa_1 \sum_{i \in \mathcal{J}} E\xi_i^2.$$

Similarly, by (3.21) with $a = \kappa_1$,

$$ET^4 \leq 3\sigma^4 + 5.5 \sum_{i \in \mathcal{J}} \left\{ \kappa_1^3 E\xi_i^4 + \kappa_1^{-1} \kappa_1^3 \sum_{j \in A_i} E|\xi_j|^4 \right.$$

$$\left. + \kappa_1^{-1} \kappa_1^3 \sum_{j \in B_i} E|\xi_j|^4 + \kappa_1^{-1} \kappa_1^3 \sum_{j \in C_i} E|\xi_j|^4 \right\}$$

$$\leq 3\sigma^4 + 22\kappa_1^3 \sum_{i \in \mathcal{J}} E|\xi_i|^4.$$

This completes the proof of Lemma 3.1. $\square$

LEMMA 3.2. *Let $\{X_i, i \in \mathcal{J}\}$ be a random field satisfying* (LD4*). *Let $\xi_i$ be a measurable function of $X_i$ with $E\xi_i = 0$ and $E\xi_i^4 < \infty$ and let $\eta_i$ be a measurable function of $X_{A_i}$ with $E\eta_i = 0$ and $E\eta_i^4 < \infty$. Let $T = \sum_{i \in \mathcal{J}} \xi_i$ and $S = \sum_{i \in \mathcal{J}} \eta_i$. Then, for any $a > 0$,*

(3.26) $\quad E(T^2 S^2) \leq 3 ET^2 ES^2 + 4 \sum_{i \in \mathcal{J}} \{a^3 E\xi_i^4 + a^{-1} E\xi_{C_i^*}^4 + a^{-1} E\xi_{D_i^*}^4$
$+ a^{-1} E\eta_{B_i^*}^4 + a^{-1} E\eta_{C_i^*}^4 + a^{-1} E\eta_{D_i^*}^4\},$

*where*

$$\xi_{C_i^*} = \sum_{j \in C_i^*} \xi_j, \qquad \xi_{D_i^*} = \sum_{j \in D_i^*} \xi_j,$$

$$\eta_{B_i^*} = \sum_{j \in B_i^*} \eta_j, \qquad \eta_{C_i^*} = \sum_{j \in C_i^*} \eta_j, \qquad \eta_{D_i^*} = \sum_{j \in D_i^*} \eta_j.$$

*In particular, we have*

$$(3.27) \qquad E(T^2 S^2) \leq 3ET^2 ES^2 + 12\kappa_2^3 \sum_{i \in \mathcal{J}} \{E|\xi_i|^4 + E|\eta_i|^4\},$$

*where* $\kappa_2 = \max_{i \in \mathcal{J}} \max(|D_i^*|, |\{j \in \mathcal{J} : i \in D_j^*\}|).$

PROOF. By (LD4$^*$), the following pairs of random variables are independent: (i) $\xi_i$ and $(T - \xi_{A_i})(S - \eta_{B_i^*})^2$; (ii) $\xi_i \xi_{A_i}$ and $(S - \eta_{B_i^*})^2$; (iii) $\xi_i \eta_{B_i^*}$ and $(T - \xi_{C_i^*})(S - \eta_{C_i^*})$. Thus we have

$$E(T^2 S^2) = \sum_{i \in \mathcal{J}} E(\xi_i T S^2)$$

$$= \sum_{i \in \mathcal{J}} E(\xi_i \{TS^2 - T(S - \eta_{B_i^*})^2 + \xi_{A_i}(S - \eta_{B_i^*})^2\})$$

$$= \sum_{i \in \mathcal{J}} 2E(\xi_i \eta_{B_i^*} TS) - \sum_{i \in \mathcal{J}} E(\xi_i \eta_{B_i^*}^2 T) + \sum_{i \in \mathcal{J}} E(\xi_i \xi_{A_i}) E(S - \eta_{B_i^*})^2$$

$$= 2\sum_{i \in \mathcal{J}} E\{\xi_i \eta_{B_i^*}(TS - (T - \xi_{C_i^*})(S - \eta_{C_i^*}))\}$$

$$+ 2\sum_{i \in \mathcal{J}} E(\xi_i \eta_{B_i^*}) E(T - \xi_{C_i^*})(S - \eta_{C_i^*})$$

$$- \sum_{i \in \mathcal{J}} E(\xi_i \eta_{B_i^*}^2 \xi_{C_i^*}) + \sum_{i \in \mathcal{J}} E(\xi_i \xi_{A_i})\{E(S^2) - 2E\eta_{B_i^*} S + E\eta_{B_i^*}^2\}$$

$$= 2\sum_{i \in \mathcal{J}} E\{\xi_i \eta_{B_i^*}(\xi_{C_i^*} S + \eta_{C_i^*} T - \xi_{C_i^*} \eta_{C_i^*})\}$$

$$+ 2\sum_{i \in \mathcal{J}} E(\xi_i \eta_{B_i^*})\{E(TS) - E(\xi_{C_i^*} S) - E(\eta_{C_i^*} T) + E(\xi_{C_i^*} \eta_{C_i^*})\}$$

$$- \sum_{i \in \mathcal{J}} E(\xi_i \eta_{B_i^*}^2 \xi_{C_i^*}) + ET^2 E(S^2)$$

$$- 2\sum_{i \in \mathcal{J}} E(\xi_i \xi_{A_i}) E\eta_{B_i^*} \eta_{C_i^*} + \sum_{i \in \mathcal{J}} E(\xi_i \xi_{A_i}) E\eta_{B_i^*}^2$$

$$= 2(E(TS))^2 + ET^2 ES^2 + 2\sum_{i \in \mathcal{J}} E\{\xi_i \eta_{B_i^*} \xi_{C_i^*} \eta_{D_i^*}\}$$

$$+ 2 \sum_{i \in \mathcal{J}} E\{\xi_i \eta_{B_i^*} \eta_{C_i^*} \xi_{D_i^*}\} - 2 \sum_{i \in \mathcal{J}} E\{\xi_i \eta_{B_i^*} \xi_{C_i^*} \eta_{C_i^*}\}$$

$$+ 2 \sum_{i \in \mathcal{J}} E(\xi_i \eta_{B_i^*})\{-E(\xi_{C_i^*} \eta_{D_i^*}) - E(\eta_{C_i^*} \xi_{D_i^*}) + E(\xi_{C_i^*} \eta_{C_i^*})\}$$

$$- \sum_{i \in \mathcal{J}} E(\xi_i \eta_{B_i^*}^2 \xi_{C_i^*}) - 2 \sum_{i \in \mathcal{J}} E(\xi_i \xi_{A_i}) E \eta_{B_i^*} \eta_{C_i^*} + \sum_{i \in \mathcal{J}} E(\xi_i \xi_{A_i}) E \eta_{B_i^*}^2$$

$$\leq 3 E T^2 E S^2 + 4 \sum_{i \in \mathcal{J}} \{a^3 E \xi_i^4 + a^{-1} E \xi_{C_i^*}^4 + a^{-1} E \xi_{D_i^*}^4$$

$$+ a^{-1} E \eta_{B_i^*}^4 + a^{-1} E \eta_{C_i^*}^4 + a^{-1} E \eta_{D_i^*}^4\},$$

where (3.25) was used to obtain the last inequality. By the $C_r$ inequality, (3.27) follows directly from (3.26). □

We are now ready to state and prove a nonuniform concentration inequality.

PROPOSITION 3.3. *Let $\{X_i, i \in \mathcal{J}\}$ be a random field satisfying* (LD4$^*$) *and put $\kappa = \max_{i \in \mathcal{J}} \max(|D_i^*|, |\{j \in \mathcal{J} : i \in D_j^*\}|)$. Let $\xi_i$ be a measurable function of $X_i$ satisfying $E\xi_i = 0$ and $|\xi_i| \leq 1/(4\kappa)$. Define*

$$\xi_{A_j} = \sum_{l \in A_j} \xi_l, \qquad \xi_{B_j} = \sum_{l \in B_j} \xi_l, \qquad T = \sum_{j \in \mathcal{J}} \xi_j.$$

*Assume that $1/2 \leq ET^2 \leq 2$ and let $\zeta = \zeta_i = (\xi_i, \xi_{A_i}, \xi_{B_i})$. Then, for Borel measurable functions $a_\zeta$ and $b_\zeta$ of $\zeta$ such that $b_\zeta \geq a_\zeta \geq 0$,*

(3.28) $\quad E^\zeta\{(1+T)^3 I(a_\zeta \leq T \leq b_\zeta)\} \leq C(b_\zeta - a_\zeta + \alpha) \quad$ a.s.,

*where $\alpha = 16\kappa^2 \sum_{j \in \mathcal{J}} E|\xi_j|^3$. In particular, if $a_\zeta \geq x > -1$, then*

(3.29) $\quad P^\zeta(a_\zeta \leq T \leq b_\zeta) \leq C(1+x)^{-3}(b_\zeta - a_\zeta + \alpha) \quad$ a.s.

PROOF. Let $C_i^{*c} = \mathcal{J} - C_i^*$ and let

$$\hat{K}_{j,\xi}(t) = \xi_j\{I(-\xi_{A_j} \leq t < 0) - I(0 \leq t \leq -\xi_{A_j})\},$$

$$\hat{M}_\xi(t) = \sum_{j \in C_i^{*c}} \hat{K}_{j,\xi}(t), \qquad M_\xi(t) = E\hat{M}_\xi(t), \qquad T_i = \sum_{j \in C_i^{*c}} \xi_j = T - \xi_{C_i^*}.$$

Since $|\xi_j| \leq 1/(4\kappa)$ and $1/2 \leq ET^2 \leq 2$, we have

(3.30) $\quad |\xi_{C_j^*}| \leq \tfrac{1}{4}, \qquad E|T - \xi_{C_j^*}|^2 \leq 6,$

and by (3.23),

(3.31) $\quad E|T - \xi_{C_j^*}|^4 \leq 108 + 22\kappa^2 \sum_{j \in \mathcal{J}} E|\xi_j|^3 \leq 108 + 2\alpha.$

Consider two cases.

CASE I ($\alpha > 1$). By (3.23),

$$E^\zeta\{(1+T)^3 I(a_\zeta \leq T \leq b_\zeta)\}$$
$$\leq (E^\zeta |1+T|^4)^{3/4} \leq \{E^\zeta (2+|T-\xi_{C_i^*}|)^4\}^{3/4}$$
$$\leq E(2+|T-\xi_{C_i^*}|)^4$$
$$\leq 8(16 + E|T-\xi_{C_i^*}|^4)$$
$$\leq 8\left(16 + 3(E|T-\xi_{C_i^*}|^2)^2 + 22\kappa^3 \sum_{j \in \mathcal{J}} E|\xi_j|^4\right)$$
$$\leq 2^{16}\left(1 + \kappa^2 \sum_{j \in \mathcal{J}} E|\xi_j|^3\right) \leq 2^{17}\alpha.$$

This proves (3.28).

CASE II ($0 < \alpha < 1$). Define

$$h_\zeta(w) = \begin{cases} 0, & \text{for } w \leq a_\zeta - \alpha, \\ \dfrac{1}{2\alpha}(w - a_\zeta + \alpha)^2, & \text{for } a_\zeta - \alpha < w \leq a_\zeta, \\ w - a_\zeta + \dfrac{\alpha}{2}, & \text{for } a_\zeta < w \leq b_\zeta, \\ -\dfrac{1}{2\alpha}(w - b_\zeta - \alpha)^2 + b_\zeta - a_\zeta + \alpha, & \text{for } b_\zeta < w \leq b_\zeta + \alpha, \\ b_\zeta - a_\zeta + \alpha, & \text{for } w > b_\zeta + \alpha, \end{cases}$$

and $f_\zeta(w) = (1+w)^3 h_\zeta(w)$. Clearly, $h'_\zeta$ is a continuous function given by

$$(3.32) \quad h'_\zeta(w) = \begin{cases} 1, & \text{for } a_\zeta \leq w \leq b_\zeta, \\ 0, & \text{for } w \leq a_\zeta - \alpha \text{ or } w \geq b_\zeta + \alpha, \\ \text{linear}, & \text{for } a_\zeta - \alpha \leq w \leq a_\zeta \text{ or } b_\zeta \leq w \leq b_\zeta + \alpha, \end{cases}$$

and $0 \leq h_\zeta(w) \leq b_\zeta - a_\zeta + \alpha$. With this $f_\zeta$, and by the fact that for every $j \in C_i^{*c}$, $\xi_j$ and $(\zeta, T - \xi_{A_j})$ are independent,

$$\begin{aligned}
E^\zeta\{T_i f_\zeta(T)\} &= \sum_{j \in C_i^{*c}} E^\zeta(\xi_j\{f_\zeta(T) - f_\zeta(T - \xi_{A_j})\}) \\
&= \sum_{j \in C_i^{*c}} E^\zeta(\xi_j\{(1+T)^3 - (1+T-\xi_{A_j})^3\}h_\zeta(T - \xi_{A_j})) \\
&\quad + \sum_{j \in C_i^{*c}} E^\zeta\{\xi_j(1+T)^3(h_\zeta(T) - h_\zeta(T - \xi_{A_j}))\} \\
&:= G_1 + G_2.
\end{aligned}$$
(3.33)

Write $G_1 = 3G_{1,1} - 3G_{1,2} + G_{1,3}$, where

$$G_{1,1} = \sum_{j \in C_i^{*c}} E^\zeta(\xi_j \xi_{A_j}(1+T)^2 h_\zeta(T - \xi_{A_j})),$$

$$G_{1,2} = \sum_{j \in C_i^{*c}} E^\zeta(\xi_j \xi_{A_j}^2 (1+T) h_\zeta(T - \xi_{A_j})),$$

$$G_{1,3} = \sum_{j \in C_i^{*c}} E^\zeta(\xi_j \xi_{A_j}^3 h_\zeta(T - \xi_{A_j})).$$

Then $(\xi_j, \xi_{A_j}, T - \xi_{C_i^*})$ and $\zeta$ are independent for each $j \in C_i^{*c}$. Hence, by (3.30),

$$|G_{1,2}| \leq (b_\zeta - a_\zeta + \alpha) \sum_{j \in C_i^{*c}} E^\zeta\{|\xi_j|\xi_{A_j}^2(1+|T|)\}$$

$$\leq (b_\zeta - a_\zeta + \alpha) \sum_{j \in C_i^{*c}} E\{|\xi_j|\xi_{A_j}^2(3+|T - \xi_{C_i^*}|)\}$$

$$\leq C(b_\zeta - a_\zeta + \alpha) \sum_{j \in C_i^{*c}} E\{|\xi_j|\xi_{A_j}^2(1+|T - \xi_{C_j^*}|)\}$$

$$\leq C(b_\zeta - a_\zeta + \alpha) \sum_{j \in C_i^{*c}} E\{|\xi_j|\xi_{A_j}^2\} E(1+|T - \xi_{C_j^*}|)$$

$$\leq C(b_\zeta - a_\zeta + \alpha) \sum_{j \in \mathcal{J}} E\{|\xi_j|\xi_{A_j}^2\}$$

$$\leq C(b_\zeta - a_\zeta + \alpha)\kappa^2 \sum_{j \in \mathcal{J}} E|\xi_j|^3$$

$$\leq C(b_\zeta - a_\zeta + \alpha).$$

Similarly, by (3.30), $|G_{1,3}| \leq C(b_\zeta - a_\zeta + \alpha)$. To bound $G_{1,1}$ write

$$G_{1,1} = \sum_{j \in C_i^{*c}} E(\xi_j \xi_{A_j}) E^\zeta((1 + T - \xi_{C_j^*})^2 h_\zeta(T - \xi_{C_j^*}))$$

$$+ \sum_{j \in C_i^{*c}} E^\zeta(\xi_j \xi_{A_j}\{(1+T)^2 h_\zeta(T - \xi_{A_j}) - (1 + T - \xi_{C_j^*})^2 h_\zeta(T - \xi_{C_j^*})\})$$

$$= E\left(\sum_{j \in C_i^{*c}} \xi_j \xi_{A_j}\right) E^\zeta((1+T)^2 h_\zeta(T))$$

$$+ \sum_{j \in C_i^{*c}} E(\xi_j \xi_{A_j}) E^\zeta((1 + T - \xi_{C_j^*})^2 h_\zeta(T - \xi_{C_j^*}) - (1+T)^2 h_\zeta(T))$$

$$+ \sum_{j \in C_i^{*c}} E^\zeta(\xi_j \xi_{A_j}\{(1+T)^2 h_\zeta(T - \xi_{A_j}) - (1 + T - \xi_{C_j^*})^2 h_\zeta(T - \xi_{C_j^*})\})$$

$$:= G_{1,1,1} + G_{1,1,2} + G_{1,1,3}.$$

First we have

$$|G_{1,1,1}| = E(T - \xi_{C_i^*})^2 E^\zeta((1+T)^2 h_\zeta(T))$$
$$\leq C(b_\zeta - a_\zeta + \alpha) E^\zeta (1+T)^2$$
$$\leq C(b_\zeta - a_\zeta + \alpha)(1 + E(1 + T - \xi_{C_i^*})^2)$$
$$\leq C(b_\zeta - a_\zeta + \alpha).$$

Next, as in bounding $G_{1,2}$, we obtain

$$|G_{1,1,3}| \leq \left| \sum_{j \in C_i^{*c}} E^\zeta(\xi_j \xi_{A_j}(1+T)^2 \{h_\zeta(T - \xi_{A_j}) - h_\zeta(T - \xi_{C_j^*})\}) \right|$$
$$+ \left| \sum_{j \in C_i^{*c}} E^\zeta(\xi_j \xi_{A_j} \{(1+T)^2 - (1 + T - \xi_{C_j^*})^2\} h_\zeta(T - \xi_{C_j^*})) \right|$$
$$\leq \sum_{j \in C_i^{*c}} E^\zeta(|\xi_j \xi_{A_j}|(1+T)^2 |\xi_{A_j} - \xi_{C_j^*}|)$$
$$+ C(b_\zeta - a_\zeta + \alpha) \sum_{j \in C_i^{*c}} E^\zeta(|\xi_j \xi_{A_j} \xi_{C_j^*}|(1 + |T|))$$
$$\leq C \sum_{j \in \mathcal{J}} (E|\xi_j| \xi_{A_j}^2 + E|\xi_j \xi_{A_j} \xi_{C_j^*}|) + C(b_\zeta - a_\zeta + \alpha) \sum_{j \in \mathcal{J}} E|\xi_j \xi_{A_j} \xi_{C_j^*}|$$
$$\leq C(b_\zeta - a_\zeta + \alpha).$$

Finally, in a similar way, $|G_{1,1,2}| \leq C(b_\zeta - a_\zeta + \alpha)$. Combining the above inequalities yields

(3.34) $$|G_1| \leq C(b_\zeta - a_\zeta + \alpha).$$

Now we bound $G_2$. Using the definition of $\hat{K}_{j,\xi}(t)$, we write

$$G_2 = \sum_{j \in C_i^{*c}} E^\zeta \left( \xi_j (1+T)^3 \int_{-\xi_{A_j}}^0 h'_\zeta(T+t) \, dt \right)$$
$$= \sum_{j \in C_i^{*c}} E^\zeta \left( (1+T)^3 \int_{|t| \leq 1} h'_\zeta(T+t) \hat{K}_{j,\xi}(t) \, dt \right)$$
$$= E^\zeta \left( (1+T)^3 \int_{|t| \leq 1} h'_\zeta(T+t) \hat{M}_\xi(t) \, dt \right)$$
$$= E^\zeta \left( (1+T)^3 \int_{|t| \leq 1} h'_\zeta(T) M_\xi(t) \, dt \right)$$

$$+ E^\zeta \left( (1+T)^3 \int_{|t|\leq 1} (h'_\zeta(T+t) - h'_\zeta(T)) M_\xi(t) \, dt \right)$$

$$+ E^\zeta \left( (1+T)^3 \int_{|t|\leq 1} h'_\zeta(T+t)(\hat{M}_\xi(t) - M_\xi(t)) \, dt \right)$$

$$:= G_{2,1} + G_{2,2} + G_{2,3}.$$

Note that

$$ET_i^2 = E(T - \xi_{C_i^*})^2 \geq \tfrac{1}{2} ET^2 - E\xi_{C_i^*}^2 \geq \tfrac{1}{4} - \tfrac{1}{8} = \tfrac{1}{8}$$

and hence

(3.35)    $G_{2,1} = ET_i^2 E^\zeta \{(1+T)^3 h'_\zeta(T)\} \geq \tfrac{1}{8} E^\zeta \{(1+T)^3 I(a_\zeta \leq T \leq b_\zeta)\}.$

By the Cauchy inequality, (3.26) and (3.31),

$$|G_{2,3}| \leq E^\zeta \left( (1+T)^4 \int_{|t|\leq 1} [h'_\zeta(T+t)]^2 \, dt \right)$$

$$+ E^\zeta \left( \int_{|t|\leq 1} (1+T)^2 (\hat{M}_\xi(t) - M_\xi(t))^2 \, dt \right)$$

$$\leq (b_\zeta - a_\zeta + 2\alpha) E^\zeta (1+T)^4$$

$$+ CE^\zeta \left\{ (1 + (T - \xi_{C_i^*})^2) \int_{|t|\leq 1} (\hat{M}_\xi(t) - M_\xi(t))^2 \, dt \right\}$$

$$\leq C(b_\zeta - a_\zeta + \alpha)(1 + E(T - \xi_{C_i^*})^4) + G_{2,3,1}$$

$$\leq C(b_\zeta - a_\zeta + \alpha) + G_{2,3,1},$$

where

$$G_{2,3,1} = \int_{|t|\leq 1} E\{(1 + (T - \xi_{C_i^*})^2)(\hat{M}_\xi(t) - M_\xi(t))^2\} \, dt.$$

By (3.22) and (3.27) we have

$$E\{(1 + (T - \xi_{C_i^*})^2)(\hat{M}_\xi(t) - M_\xi(t))^2\}$$

$$\leq C\kappa \sum_{j \in \mathcal{J}} E|\hat{K}_{j,\xi}(t)|^2 + C\kappa^3 \sum_{j \in \mathcal{J}} E|\xi_j|^4 + C\kappa^3 \sum_{j \in \mathcal{J}} E|\hat{K}_{j,\xi}(t)|^4$$

and hence

$$G_{2,3,1} \leq C\kappa \sum_{j \in \mathcal{J}} E|\xi_j|^2 |\xi_{A_j}| + C\kappa^3 \sum_{j \in \mathcal{J}} E|\xi_j|^4 + C\kappa^3 \sum_{j \in \mathcal{J}} E|\xi_j|^4 |\xi_{A_j}|$$

$$\leq C\kappa \sum_{j \in \mathcal{J}} E|\xi_j|^2 |\xi_{A_j}| + C\kappa^2 \sum_{j \in \mathcal{J}} E|\xi_j|^3$$

$$\leq C\kappa \sum_{j\in\mathcal{J}}\sum_{l\in A_j} E|\xi_j|^2|\xi_l| + C\kappa^2 \sum_{j\in\mathcal{J}} E|\xi_j|^3$$

$$\leq C\kappa^2 \sum_{j\in\mathcal{J}} E|\xi_j|^3 \leq C\alpha.$$

To bound $G_{2,2}$, define

$$L_\zeta(\alpha) = \lim_{k\to\infty} \sup_{x\geq 0, x\in Q} E^\zeta\{(1+T)^3 I(x-1/k \leq T \leq x+1/k+\alpha)\},$$

where $Q$ is the set of rational numbers and $E^\zeta$ is regarded as a regular conditional expectation given $\zeta$. Then, for $a_\zeta > 1$, so that $a_\zeta - \alpha > 0$, we have

$$G_{2,2} = E^\zeta\left((1+T)^3 \int_0^1 \int_0^t h_\zeta''(T+s)\,ds\,M_\xi(t)\,dt\right)$$
$$+ E^\zeta\left((1+T)^3 \int_{-1}^0 \int_t^0 h_\zeta''(T+s)\,ds\,M_\xi(t)\,dt\right)$$
$$= \alpha^{-1} \int_0^1 \int_0^t E^\zeta\{(1+T)^3(I(a_\zeta - \alpha \leq T+s \leq a_\zeta)$$
$$- I(b_\zeta \leq T+s \leq b_\zeta+\alpha))\}\,ds\,M_\xi(t)\,dt$$
$$+ \alpha^{-1} \int_{-1}^0 \int_t^0 E^\zeta\{(1+T)^3(I(a_\zeta - \alpha \leq T+s \leq a_\zeta)$$
$$- I(b_\zeta \leq T+s \leq b_\zeta+\alpha))\}\,ds\,M_\xi(t)\,dt$$

and

$$|G_{2,2}| \leq \alpha^{-1} \int_0^1 \int_0^t L_\zeta(\alpha)\,ds\,|M_\xi(t)|\,dt$$
$$+ \alpha^{-1} \int_{-1}^0 \int_t^0 L_\zeta(\alpha)\,ds\,|M_\xi(t)|\,dt$$
$$= \alpha^{-1} L_\zeta(\alpha) \int_{|t|\leq 1} |tM_\xi(t)|\,dt$$
$$\leq \tfrac{1}{2}\alpha^{-1} L_\zeta(\alpha) \sum_{j\in\mathcal{J}} E|\xi_j \xi_{A_j}^2|$$
$$\leq \tfrac{1}{2}\kappa^2 \alpha^{-1} L_\zeta(\alpha) \sum_{j\in\mathcal{J}} E|\xi_j|^3$$
$$\leq \tfrac{1}{16} L_\zeta(\alpha).$$

(3.36)

Therefore,

(3.37)
$$G_2 \geq \tfrac{1}{8} E^\zeta\{(1+T)^3 I(a_\zeta \leq T \leq b_\zeta)\}$$
$$- C(b_\zeta - a_\zeta + \alpha) - \tfrac{1}{16} L_\zeta(\alpha).$$

Now by (3.31)

$$
(3.38) \quad \begin{aligned} E^\zeta \{T_i f_\zeta(T)\} &\leq (b_\zeta - a_\zeta + \alpha) E^\zeta |T_i (1+T)^3| \\ &\leq C(b_\zeta - a_\zeta + \alpha) E(|T_i|(1+|T_i|^3)) \\ &\leq C(b_\zeta - a_\zeta + \alpha). \end{aligned}
$$

So combining (3.33), (3.34), (3.37) and (3.38), we have for $a_\zeta > 1$,

$$
(3.39) \quad E^\zeta\{(1+T)^3 I(a_\zeta \leq T \leq b_\zeta)\} \leq C(b_\zeta - a_\zeta + \alpha) + \tfrac{1}{2} L_\zeta(\alpha).
$$

For $0 < a_\zeta \leq 1$, it suffices to consider $b_\zeta - a_\zeta \leq 1$. Applying Proposition 3.2 to $\{\xi_i, i \in \mathcal{J}\}$, we obtain

$$
(3.40) \quad \begin{aligned} &E^\zeta\{(1+T)^3 I(a_\zeta < T < b_\zeta + \alpha)\} \\ &\qquad \leq C P^\zeta(a_\zeta < T < b_\zeta + \alpha) \leq C(b_\zeta - a_\zeta + \alpha). \end{aligned}
$$

Now take $a_\zeta = x - 1/k$ and $b_\zeta = x + 1/k + \alpha$. By taking the supremum over $x \in Q$ and then letting $k \to \infty$, (3.39) and (3.40) imply

$$
L_\zeta(\alpha) \leq C\alpha + \tfrac{1}{2} L_\zeta(\alpha).
$$

Hence

$$
(3.41) \quad L_\zeta(\alpha) \leq C\alpha.
$$

This together with (3.39) and (3.40) proves (3.28) and hence Proposition 3.3. □

**4. Proofs of Theorems 2.1–2.4.** We first derive a Stein identity for $W$. Let $f$ be a bounded absolutely continuous function. Then

$$
(4.1) \quad \begin{aligned} E\{W f(W)\} &= \sum_{i \in \mathcal{J}} E\{X_i(f(W) - f(W - Y_i))\} \\ &= \sum_{i \in \mathcal{J}} E\left\{ X_i \int_{-Y_i}^{0} f'(W+t)\,dt \right\} \\ &= \sum_{i \in \mathcal{J}} E\left\{ \int_{-\infty}^{\infty} f'(W+t) \hat{K}_i(t)\,dt \right\} \\ &= E \int_{-\infty}^{\infty} f'(W+t) \hat{K}(t)\,dt \end{aligned}
$$

and hence by the fact that $\int_{-\infty}^{\infty} K(t)\,dt = EW^2 = 1$,

$$
\begin{aligned}
Ef'(W) &- EWf(W) \\
&= E\int_{-\infty}^{\infty} f'(W)K(t)\,dt - E\int_{-\infty}^{\infty} f'(W+t)\hat{K}(t)\,dt \\
&= E\int_{-\infty}^{\infty} f'(W)(K(t) - \hat{K}(t))\,dt \\
&\quad + E\int_{|t|>1} (f'(W) - f'(W+t))\hat{K}(t)\,dt \\
&\quad + E\int_{|t|\leq 1} (f'(W) - f'(W+t))(\hat{K}(t) - K(t))\,dt \\
&\quad + E\int_{|t|\leq 1} (f'(W) - f'(W+t))K(t)\,dt \\
&:= R_1 + R_2 + R_3 + R_4.
\end{aligned}
\tag{4.2}
$$

Now choose $f$ to be $f_{z,\alpha}$, the unique bounded solution of the differential equation

$$f'(w) - wf(w) = h_{z,\alpha}(w) - Nh_{z,\alpha}, \tag{4.3}$$

where $\alpha > 0$ is to be determined later,

$$h_{z,\alpha}(w) = \begin{cases} 1, & \text{for } w \leq z, \\ 1 + (z-w)/\alpha, & \text{for } z \leq w \leq z+\alpha, \\ 0, & \text{for } w \geq z+\alpha, \end{cases} \tag{4.4}$$

and $Nh_{z\alpha} = (2\pi)^{-1/2} \int_{-\infty}^{\infty} h_{z,\alpha}(x) e^{-x^2/2}\,dx$. The solution of (4.3) is given by

$$
\begin{aligned}
f_{z,\alpha}(w) &= e^{w^2/2} \int_{-\infty}^{w} [h_{z,\alpha}(x) - Nh_{z,\alpha}] e^{-x^2/2}\,dx \\
&= -e^{w^2/2} \int_{w}^{\infty} [h_{z,\alpha}(x) - Nh_{z,\alpha}] e^{-x^2/2}\,dx.
\end{aligned}
$$

From Lemma 2 and arguments on pages 23 and 24 in Stein (1986), we have, for all $w$ and $v$,

$$0 \leq f_{z,\alpha}(w) \leq 1, \tag{4.5}$$

$$|f'_{z,\alpha}(w)| \leq 1, \qquad |f'_{z,\alpha}(w) - f'_{z,\alpha}(v)| \leq 1 \tag{4.6}$$

and

$$
\begin{aligned}
&|f'_{z,\alpha}(w+s) - f'_{z,\alpha}(w+t)| \\
&\qquad \leq (|w|+1)\min(|s|+|t|,1) + \alpha^{-1}\left|\int_s^t I(z \leq w+u \leq z+\alpha)\,du\right| \quad (4.7) \\
&\qquad \leq (|w|+1)\min(|s|+|t|,1) + I(z - s \vee t \leq w \leq z - s \wedge t + \alpha). \quad (4.8)
\end{aligned}
$$

PROOF OF THEOREM 2.1. By (4.6),

$$|R_1| = \left|Ef'(W)\sum_{i\in\mathcal{J}}(X_iY_i - EX_iY_i)\right| \leq r_1 \tag{4.9}$$

and

$$|R_2| \leq \sum_{i\in\mathcal{J}} E|X_iY_i|I(|Y_i| > 1) = r_2. \tag{4.10}$$

By (4.8),

$$\begin{aligned}|R_3| &\leq E\int_{|t|\leq 1}(|W|+1)|t||\hat{K}(t) - K(t)|\,dt \\ &\quad + E\int_0^1 I(z - t \leq W \leq z + \alpha)|\hat{K}(t) - K(t)|\,dt \\ &\quad + E\int_{-1}^0 I(z \leq W \leq z - t + \alpha)|\hat{K}(t) - K(t)|\,dt \\ &\leq r_3 + r_4 + R_{3,1} + R_{3,2},\end{aligned} \tag{4.11}$$

where

$$R_{3,1} = E\int_0^1 I(z - t \leq W \leq z + \alpha)|\hat{K}(t) - K(t)|\,dt,$$

$$R_{3,2} = E\int_{-1}^0 I(z \leq W \leq z - t + \alpha)|\hat{K}(t) - K(t)|\,dt.$$

Let $\delta = 0.625\alpha + 4r_2 + 2.125r_3 + 4r_5$. Then by Proposition 3.1

$$P(z - t \leq W \leq z + \alpha) \leq \delta + 0.625t$$

for $t > 0$. Hence by the Cauchy inequality,

$$R_{3,1} \leq E\bigg\{\int_0^1 (0.5\alpha(\delta + 0.625t)^{-1}I(z - t \leq W \leq z + \alpha)$$

$$+ 0.5\alpha^{-1}(\delta + 0.625t)|\hat{K}(t) - K(t)|^2)\,dt\bigg\}$$

$$\leq 0.5\alpha + 0.5\alpha^{-1}\delta\int_0^1 \mathrm{Var}(\hat{K}(t))\,dt + 0.32\alpha^{-1}\int_0^1 t\,\mathrm{Var}(\hat{K}(t))\,dt.$$

A similar inequality holds for $R_{3,2}$. Thus we arrive at

$$R_3 \leq \alpha + 0.5\alpha^{-1}\delta r_5 + 0.32\alpha^{-1}r_6^2 + r_3 + r_4. \tag{4.12}$$

By (4.7) and Proposition 3.1 again, we have

$$\begin{aligned}|R_4| &\leq E\int_{|t|\leq 1}(|W|+1)|tK(t)|\,dt \\ &\quad + \alpha^{-1}\int_{|t|\leq 1}\left|\int_0^t P(z \leq W + u \leq z + \alpha)\,du\right||K(t)|\,dt \\ &\leq 2r_3 + \alpha^{-1}\int_{|t|\leq 1} t\delta|K(t)|\,dt \leq 2r_3 + \alpha^{-1}\delta r_3.\end{aligned} \tag{4.13}$$

NORMAL APPROXIMATION UNDER LOCAL DEPENDENCE 27Combining the above inequalities yields

$$|Eh_{z,\alpha}(W) - Nh_{z,\alpha}|$$
$$\leq r_1 + r_2 + 3r_3 + r_4 + \alpha + \alpha^{-1}\{\delta(0.5r_5 + r_3) + 0.32r_6^2\}$$
$$\leq r_1 + r_2 + 3.625r_3 + r_4 + 0.32r_5 + \alpha$$
$$+ \alpha^{-1}\{(4r_2 + 2.125r_3 + 4r_5)(0.5r_5 + r_3) + 0.32r_6^2\}.$$

Using the fact that $Eh_{z-\alpha,\alpha}(W) \leq P(W \leq z) \leq Eh_{z,\alpha}(W)$ and that $|\Phi(z+\alpha) - \Phi(z)| \leq (2\pi)^{-1/2}\alpha$, we have

(4.14) $$\sup_z |P(W \leq z) - \Phi(z)| \leq \sup_z |Eh_{z,\alpha}(W) - Nh_{z,\alpha}| + 0.5\alpha.$$

Letting $\alpha = \sqrt{2/3}((4r_2 + 2.125r_3 + 4r_5)(0.5r_5 + r_3) + 0.32r_6^2)^{1/2}$ yields

(4.15)
$$\sup_z |P(W \leq z) - \Phi(z)|$$
$$\leq r_1 + r_2 + 3.625r_3 + r_4 + 0.32r_5$$
$$+ \sqrt{6}((4r_2 + 2.125r_3 + 4r_5)(0.5r_5 + r_3) + 0.32r_6^2)^{1/2}$$
$$\leq r_1 + r_2 + 3.625r_3 + r_4 + 0.32r_5 + 1.5r_6$$
$$+ 0.5(1.5(4r_2 + 2.125r_3 + 4r_5) + 4(0.5r_5 + r_3))$$
$$\leq r_1 + 4r_2 + 8r_3 + r_4 + 4.5r_5 + 1.5r_6.$$

This proves Theorem 2.1. □

PROOF OF THEOREM 2.2. Let $p_3 = 3 \wedge p$. Using the following well-known inequality: $\forall x_i \geq 0, \alpha_i \geq 0$ with $\sum \alpha_i = 1$,

(4.16) $$\prod x_i^{\alpha_i} \leq \sum \alpha_i x_i,$$

we have

$$r_2 \leq \sum_{i \in \mathcal{J}} E|X_i||Y_i|^{p_3-1}$$
$$\leq \sum_{i \in \mathcal{J}} \left\{ \frac{1}{p_3} E|X_i|^{p_3} + \frac{p_3-1}{p_3} E|Y_i|^{p_3} \right\}$$

and

$$r_3 \leq \sum_{i \in \mathcal{J}} E|X_i||Y_i|^{p_3-1}$$
$$\leq \sum_{i \in \mathcal{J}} \left\{ \frac{1}{p_3} E|X_i|^{p_3} + \frac{p_3-1}{p_3} E|Y_i|^{p_3} \right\}.$$

Similarly we have

$$r_5 \leq \sum_{i,j \in \mathcal{J}, B_i B_j \neq \varnothing} \{E|X_i X_j||Y_i|^{p_3-2} + E|X_i X_j^*||Y_i|^{p_3-2}\}$$

$$\leq 2 \sum_{i,j \in \mathcal{J}, B_i B_j \neq \varnothing} \left\{ \frac{1}{p_3} E|X_i|^{p_3} + \frac{1}{p_3} E|X_j|^{p_3} + \frac{p_3-2}{p_3} E|Y_i|^{p_3} \right\}$$

$$\leq 2\kappa \sum_{i \in \mathcal{J}} \left\{ \frac{2}{p_3} E|X_i|^{p_3} + \frac{p_3-2}{p_3} E|Y_i|^{p_3} \right\}$$

and

$$r_6^2 \leq \frac{1}{2} \sum_{i,j \in \mathcal{J}, B_i B_j \neq \varnothing} \{E|X_i X_j||Y_i|^{p-2} + E|X_i X_j^*||Y_i|^{p-2}\}$$

$$\leq \kappa \sum_{i \in \mathcal{J}} \left\{ \frac{2}{p} E|X_i|^p + \frac{p-2}{p} E|Y_i|^p \right\}.$$

Now we estimate $r_4$. Recall that

$$Z_i = \sum_{j \in B_i} X_j$$

and that $(X_i, Y_i)$ and $W - Z_i$ are independent.

We have

$$r_4 \leq \sum_{i \in \mathcal{J}} \{E|W - Z_i|E|X_i|(Y_i^2 \wedge 1) + E|Z_i X_i|(Y_i^2 \wedge 1)\}$$

$$\leq \sum_{i \in \mathcal{J}} \{(1 + E|Z_i|)E|X_i|(Y_i^2 \wedge 1) + E|Z_i X_i||Y_i|^{p_3-2}\}$$

$$\leq \sum_{i \in \mathcal{J}} E|X_i||Y_i|^{p_3-1} + \sum_{i \in \mathcal{J}} \sum_{j \in B_i} E|X_j|E|X_i||Y_i|^{p_3-2}$$

$$+ \sum_{i \in \mathcal{J}} \sum_{j \in B_i} E|X_j X_i||Y_i|^{p_3-2}$$

$$\leq \sum_{i \in \mathcal{J}} \left\{ \frac{1}{p_3} E|X_i|^{p_3} + \frac{p_3-1}{p_3} E|Y_i|^{p_3} \right\}$$

$$+ 2\kappa \sum_{i \in \mathcal{J}} \left\{ \frac{2}{p_3} E|X_i|^{p_3} + \frac{p_3-2}{p_3} E|Y_i|^{p_3} \right\}.$$

To estimate $r_1$, let $\xi_i = X_i Y_i I(|X_i Y_i| \leq 1)$. We have

$$r_1 \leq E\left|\sum_{i \in \mathcal{J}}(\xi_i - E\xi_i)\right| + 2\sum_{i \in \mathcal{J}} E|X_i Y_i|I(|X_i Y_i| > 1)$$

$$\leq \operatorname{Var}\left(\sum_{i \in \mathcal{J}} \xi_i\right)^{1/2} + 2 \sum_{i \in \mathcal{J}} E|X_i Y_i|^{p_3/2}$$

$$\leq \operatorname{Var}\left(\sum_{i \in \mathcal{J}} \xi_i\right)^{1/2} + \sum_{i \in \mathcal{J}} \{E|X_i|^{p_3} + E|Y_i|^{p_3}\}.$$

Similarly to bounding $r_6$,

$$\operatorname{Var}\left(\sum_{i \in \mathcal{J}} \xi_i\right) \leq \sum_{i,j \in \mathcal{J}, B_i B_j \neq \varnothing} (E|\xi_i \xi_j| + E|\xi_i| E|\xi_j|)$$

$$\leq \sum_{i,j \in \mathcal{J}, B_i B_j \neq \varnothing} (E|X_i Y_i|^{p/4} |X_j Y_j|^{p/4} + E|X_i Y_i|^{p/4} E|X_j Y_j|^{p/4})$$

$$\leq \kappa \sum_{i \in \mathcal{J}} (E|X_i|^p + E|Y_i|^p).$$

Combining the inequalities above yields (2.5). □

PROOF OF THEOREM 2.3. The idea of the proof is similar to that of Theorem 2.1. Noting that $\hat{K}_i$ and $W - Z_i$ are independent, we rewrite (4.2) as

$$\begin{aligned}
Ef'(W) &- EWf(W) \\
&= \sum_{i \in \mathcal{J}} \left\{ E \int_{-\infty}^{\infty} f'(W) K_i(t)\, dt - E \int_{-\infty}^{\infty} f'(W+t) \hat{K}_i(t)\, dt \right\} \\
&= \sum_{i \in \mathcal{J}} E \int_{-\infty}^{\infty} (f'(W) - f'(W - Z_i + t)) K_i(t)\, dt \\
&\quad + \sum_{i \in \mathcal{J}} E \int_{-\infty}^{\infty} (f'(W - Z_i + t) - f'(W + t)) \hat{K}_i(t)\, dt \\
&= Q_1 + Q_2 + Q_3 + Q_4,
\end{aligned}$$
(4.17)

where

$$Q_1 = \sum_{i \in \mathcal{J}} E \int_{|t| \leq 1} (f'(W) - f'(W - Z_i + t)) K_i(t)\, dt,$$

$$Q_2 = \sum_{i \in \mathcal{J}} E \int_{|t| > 1} (f'(W) - f'(W - Z_i + t)) K_i(t)\, dt,$$

$$Q_3 = \sum_{i \in \mathcal{J}} E \int_{|t| > 1} (f'(W - Z_i + t) - f'(W + t)) \hat{K}_i(t)\, dt,$$

$$Q_4 = \sum_{i \in \mathcal{J}} E \int_{|t| \leq 1} (f'(W - Z_i + t) - f'(W + t)) \hat{K}_i(t)\, dt.$$

By (4.6), similarly to (4.10),

(4.18) $$|Q_2| + |Q_3| \leq 2 \sum_{i \in \mathcal{J}} E|X_i Y_i| I(|Y_i| > 1) = 2r_2.$$

To bound $Q_4$, write $Q_4 = Q_{4,1} + Q_{4,2}$, where

$$Q_{4,1} = \sum_{i \in \mathcal{J}} EI(|X_i| > 1) \int_{|t| \leq 1} (f'(W - Z_i + t) - f'(W + t)) \hat{K}_i(t) \, dt,$$

$$Q_{4,2} = \sum_{i \in \mathcal{J}} EI(|X_i| \leq 1) \int_{|t| \leq 1} (f'(W - Z_i + t) - f'(W + t)) \hat{K}_i(t) \, dt.$$

Then by (4.6),

(4.19) $$|Q_{4,1}| \leq \sum_{i \in J} E|X_i Y_i| I(|X_i| > 1) = r_7.$$

From (4.7), we obtain

$$|Q_{4,2}| \leq \sum_{i \in \mathcal{J}} EI(|X_i| \leq 1) \int_{|t| \leq 1} (|W| + |t| + 1)(|Z_i| \wedge 1) |\hat{K}_i(t)| \, dt$$
$$+ \alpha^{-1} \sum_{i \in \mathcal{J}} EI(|X_i| \leq 1)$$
$$\times \int_{|t| \leq 1} I(Z_i \geq 0)$$
$$\times \int_{-Z_i}^{0} I(z \leq W + t + u \leq z + \alpha) \, du \, |\hat{K}_i(t)| \, dt$$
$$+ \alpha^{-1} \sum_{i \in \mathcal{J}} EI(|X_i| \leq 1)$$
$$\times \int_{|t| \leq 1} I(Z_i < 0)$$
$$\times \int_{0}^{-Z_i} I(z \leq W + t + u \leq z + \alpha) \, du |\hat{K}_i(t)| \, dt$$
$$\leq \sum_{i \in \mathcal{J}} E(|W| + 1)|X_i| I(|X_i| \leq 1)(|Z_i| \wedge 1)(|Y_i| \wedge 1)$$
$$+ 0.5 \sum_{i \in \mathcal{J}} E|X_i|(|Y_i|^2 \wedge 1) + Q_{4,3}$$
$$\leq r_8 + r_9 + 0.5 r_3 + Q_{4,3},$$
(4.20)
where

$$Q_{4,3} = \alpha^{-1} \sum_{i \in \mathcal{J}} E \Big\{ I(|X_i| \leq 1)$$
$$\times \int_{|t| \leq 1} I(Z_i \geq 0)$$

$$\times \int_{-Z_i}^0 P^{\xi_i}(z \leq W + t + u \leq z + \alpha)\,du\,|\hat{K}_i(t)|\,dt\Big\}$$

$$+ \alpha^{-1}\sum_{i\in\mathcal{J}} E\Big\{I(|X_i| \leq 1)$$

$$\times \int_{|t|\leq 1} I(Z_i < 0)$$

$$\times \int_0^{-Z_i} P^{\xi_i}(z \leq W + t + u \leq z + \alpha)\,du\,|\hat{K}_i(t)|\,dt\Big\}$$

and $\xi_i = (X_i, Y_i, Z_i)$. By Proposition 3.2,

$$Q_{4,3} \leq \alpha^{-1}\sum_{i\in\mathcal{J}} E\Big\{I(|X_i| \leq 1)$$

$$\times \int_{|t|\leq 1} (0.625\sigma_i^{-1}\alpha + 4\sigma_i^{-2}r_2$$

(4.21)
$$+ 2.125\sigma_i^{-3}r_3 + 4\sigma_i^{-3}r_{10})|Z_i||\hat{K}_i(t)|\,dt\Big\}$$

$$\leq 0.625\lambda \sum_{i\in\mathcal{J}} E\{|X_i|I(|X_i|\leq 1)(|Y_i|\wedge 1)|Z_i|\}$$

$$+ \alpha^{-1}\{4\lambda^2 r_2 + 2.125\lambda^3 r_3 + 4\lambda^3 r_{10}\}r_8$$

$$= 0.625\lambda r_8 + \alpha^{-1}\{4\lambda^2 r_2 + 2.125\lambda^3 r_3 + 4\lambda^3 r_{10}\}r_8.$$

Combining (4.19)–(4.21) yields

(4.22) $$|Q_4| \leq r_7 + r_9 + 1.625\lambda r_8 + 0.5 r_3 \\ + \alpha^{-1}\{4\lambda^2 r_2 + 2.125\lambda^3 r_3 + 4\lambda^3 r_{10}\}r_8.$$

Similarly we have

(4.23)
$$|Q_1| \leq \sum_{i\in\mathcal{J}} P(|X_i| > 1)E|X_i|(|Y_i|\wedge 1) \\ + \sum_{i\in\mathcal{J}} E\{(|W|+1)(|Z_i|\wedge 1)\}E|X_i|(|Y_i|\wedge 1) + r_3 \\ + \alpha^{-1}\{0.625\lambda\alpha + 4\lambda^2 r_2 + 2.125\lambda^3 r_3 + 4\lambda^3 r_{10}\} \\ \times \Big\{0.5 r_3 + \sum_{i\in\mathcal{J}} E\{(|W|+1)(|Z_i|\wedge 1)\}E|X_i|(|Y_i|\wedge 1)\Big\} \\ = r_{11} + r_{12} + r_3 \\ + \alpha^{-1}\{0.625\lambda\alpha + 4\lambda^2 r_2 + 2.125\lambda^3 r_3 + 4\lambda^3 r_{10}\}(0.5 r_3 + r_{12}).$$

Combining (4.17), (4.18), (4.22), (4.23) and (4.14), we obtain

(4.24)
$$\sup_z |F(z) - \Phi(z)|$$
$$\leq 0.5\alpha + 2r_2 + 2\lambda r_3 + r_7 + 1.625\lambda r_8 + r_9 + r_{11} + 1.625\lambda r_{12}$$
$$+ \alpha^{-1}\{4\lambda^2 r_2 + 2.125\lambda^3 r_3 + 4\lambda^3 r_{10}\}(r_8 + 0.5r_3 + r_{12}).$$

Let

$$\alpha = (2(4\lambda^2 r_2 + 2.125\lambda^3 r_3 + 4\lambda^3 r_{10})(r_8 + 0.5r_3 + r_{12}))^{1/2}.$$

Then the right-hand side of (4.24) is

$$= 2r_2 + 2\lambda r_3 + r_7 + 1.625\lambda r_8 + r_9 + r_{11} + 1.625\lambda r_{12}$$
$$+ \{2(4\lambda^2 r_2 + 2.125\lambda^3 r_3 + 4\lambda^3 r_{10})(r_8 + 0.5r_3 + r_{12})\}^{1/2}$$
$$\leq 2r_2 + 2\lambda r_3 + r_7 + 1.625\lambda r_8 + r_9 + r_{11} + 1.625\lambda r_{12}$$
$$+ 0.5\lambda^{-3/2}(4\lambda^2 r_2 + 2.125\lambda^3 r_3 + 4\lambda^3 r_{10}) + \lambda^{3/2}(r_8 + 0.5r_3 + r_{12})$$
$$\leq 4\lambda^{3/2}(r_2 + r_3 + r_7 + r_8 + r_9 + r_{10} + r_{11} + r_{12}).$$

This proves Theorem 2.3. □

PROOF OF THEOREM 2.4. We can assume that

(4.25)
$$\kappa^{p-1} \sum_{i \in \mathcal{J}} E|X_i|^p \leq \tfrac{1}{75}.$$

Otherwise (2.9) is trivial. Let $\xi_i = \sum_{j \in N(C_i)} X_j$. Then by (4.25)

(4.26)
$$\sqrt{E\xi_i^2} \leq (E|\xi_i|^p)^{1/p} \leq \left(|N(C_i)|^{p-1} \sum_{j \in N(C_i)} E|X_j|^p\right)^{1/p}$$
$$\leq \left(\kappa^{p-1} \sum_{j \in \mathcal{J}} E|X_j|^p\right)^{1/p} \leq (\tfrac{1}{75})^{1/p} \leq (\tfrac{1}{75})^{1/3} < 0.2372$$

and

$$\sigma_i \geq \sqrt{EW^2} - \sqrt{E\xi_i^2} \geq 1 - 0.2372 = 0.7628.$$

Thus $4\lambda^{3/2} \leq 6.01$. By (4.16) and the Minkowski inequality

$$r_8 \leq \sum_{i \in \mathcal{J}} E|X_i Z_i||Y_i|^{p-2}$$
$$= \sum_{i \in \mathcal{J}} E(\kappa^{(p-1)/p}|X_i|(|Z_i|/\kappa^{1/p})(|Y_i|/\kappa^{1/p})^{p-2})$$
$$\leq \sum_{i \in \mathcal{J}} \frac{1}{p}\kappa^{p-1}E|X_i|^p + \sum_{i \in \mathcal{J}} \left(\frac{1}{p\kappa}E|Z_i|^p + \frac{(p-2)}{p\kappa}E|Y_i|^p\right)$$

$$\leq \frac{\kappa^{p-1}}{p} \sum_{i \in \mathcal{J}} E|X_i|^p + \sum_{i \in \mathcal{J}} \frac{|C_i|^{p-1}(p-1)}{p\kappa} \sum_{j \in C_i} E|X_j|^p$$

$$\leq \kappa^{p-1} \sum_{i \in \mathcal{J}} E|X_i|^p.$$

Similarly we have

$$r_2 + r_3 + r_7 + r_{11} \leq 2\kappa^{p-1} \sum_{i \in \mathcal{J}} E|X_i|^p.$$

Note that $|N(B_i)| \leq |N(C_i)|$. Following the proof for $r_6^2$ yields

$$r_{10} \leq \sum_{i,j \in \mathcal{J}, B_i B_j \neq \varnothing} \{ E(\kappa^{(p-2)/p}|X_i|\kappa^{(p-2)/p}|X_j|(|Y_i|/\kappa^{2/p})^{p-2})$$

$$+ E(\kappa^{(p-2)/p}|X_i|\kappa^{(p-2)/p}|X_j^*|(|Y_i|/\kappa^{2/p})^{p-2}) \}$$

$$\leq 2 \sum_{i,j \in \mathcal{J}, B_i B_j \neq \varnothing} \left\{ \frac{\kappa^{p-2}}{p}(E|X_i|^p + E|X_j|^p) + \frac{p-2}{p\kappa^2} E|Y_i|^p \right\}$$

$$\leq 2\kappa^{p-1} \sum_{i \in \mathcal{J}} E|X_i|^p.$$

The $r_9$ can be bounded as $r_8$ and $r_4$. By (4.26),

$$r_9 \leq \sum_{i \in \mathcal{J}} E|W - \xi_i| E|X_i|(|Y_i| \wedge 1)(|Z_i| \wedge 1) + \sum_{i \in \mathcal{J}} E|\xi_i||X_i||Y_i|^{p-2}$$

$$\leq 1.2372 \sum_{i \in \mathcal{J}} E|X_i||Y_i||Z_i|^{p-2} + \sum_{i \in \mathcal{J}} E|X_i||\xi_i||Y_i|^{p-2}$$

$$\leq 2.2372\kappa^{p-1} \sum_{i \in \mathcal{J}} E|X_i|^p.$$

Similarly,

$$r_{12} \leq 3.2372\kappa^{p-1} \sum_{i \in \mathcal{J}} E|X_i|^p.$$

Theorem 2.4 follows from (2.7) and the above inequalities. □

**5. Proof of Theorem 2.5.** The basic idea of the proof of Theorem 2.5 is similar to that of Theorem 2.3. We use the same notation as in Section 2.2 and remind the reader that $\{X_i, i \in \mathcal{J}\}$ satisfies (LD4*) and that $E|X_i|^p < \infty$ for $2 < p \leq 3$.

First we need a few preliminary lemmas. Let

(5.1) $$\tau = 1/(8\kappa), \quad \bar{X}_i = X_i I(|X_i| \leq \tau), \quad \hat{X}_i = X_i I(|X_i| > \tau),$$
$$\bar{W} = \sum_{i \in \mathcal{J}} \bar{X}_i, \quad Y_i = \sum_{j \in A_i} X_j,$$

and

$$\beta_1 = \sum_{i \in \mathcal{J}} E|X_i Y_i| I(|X_i| > \tau), \qquad \beta_2 = \sum_{i \in \mathcal{J}} E X_i^2 I(|X_i| > \tau),$$
(5.2)
$$\beta_3 = \sum_{i \in \mathcal{J}} E|X_i|^3 I(|X_i| \leq \tau).$$

Our first lemma shows that $W$ is close to $\bar{W}$.

LEMMA 5.1. *Assume $\beta_2 \leq \tau/16$. Then there exists an absolute constant $C$ such that, for $z \geq 0$,*

(5.3)
$$\begin{aligned} |P(W > z) - P(\bar{W} > z)| \\ \leq \sum_{i \in \mathcal{J}} P(|Y_i| > (z+1)/4) + C(z+1)^{-3} \kappa^3 \beta_3 \\ + 64 \kappa^2 \beta_2 \{(1+z)^{-3} + P(W > (z-1)/2)\}. \end{aligned}$$

PROOF. Observe that

(5.4)
$$\begin{aligned} P(W > z) &= P\Big(W > z, \max_{i \in \mathcal{J}} |X_i| \leq \tau\Big) + P\Big(W > z, \max_{i \in \mathcal{J}} |X_i| > \tau\Big) \\ &\leq P(\bar{W} > z) + \sum_{i \in \mathcal{J}} P(W > z, |X_i| > \tau) \\ &\leq P(\bar{W} > z) + \sum_{i \in \mathcal{J}} \{P(W - Y_i > 3(z - 1/3)/4, |X_i| > \tau) \\ &\qquad\qquad\qquad + P(Y_i > (z+1)/4, |X_i| > \tau)\} \\ &\leq P(\bar{W} > z) + \sum_{i \in \mathcal{J}} \{P(W - Y_i > 3(z - 1/3)/4) P(|X_i| > \tau) \\ &\qquad\qquad\qquad + P(Y_i > (z+1)/4)\} \\ &\leq P(\bar{W} > z) + \sum_{i \in \mathcal{J}} \{P(W > (z-1)/2) P(|X_i| > \tau) \\ &\qquad\qquad\qquad + P(-Y_i > (z+1)/4) + P(Y_i > (z+1)/4)\} \\ &\leq P(\bar{W} > z) + P(W > (z-1)/2) \tau^{-2} \sum_{i \in \mathcal{J}} E|X_i^2| I(|X_i| > \tau) \\ &\qquad + \sum_{i \in J} P(|Y_i| > (z+1)/4) \\ &= P(\bar{W} > z) + 64 \kappa^2 \beta_2 P(W > (z-1)/2) \\ &\qquad + \sum_{i \in J} P(|Y_i| > (z+1)/4). \end{aligned}$$

Similarly, noting that $|\sum_{j \in A_i} \bar{X}_j| \leq 1$, we have

(5.5) $$P(\bar{W} > z) \leq P(W > z) + \kappa^2 \beta_2 P(\bar{W} > z - 2).$$

Note that $\beta_2 \leq \tau/16$ implies $|E\bar{W}| \leq 1/16$ and by (3.22),

$$
\begin{aligned}
|\operatorname{Var}(\bar{W}) - 1| &= |\operatorname{Var}(\bar{W}) - \operatorname{Var}(W)| \\
&\leq 2(EW^2)^{1/2}\left(\operatorname{Var}\left(\sum_{i \in \mathcal{J}} \hat{X}_i\right)\right)^{1/2} + \operatorname{Var}\left(\sum_{i \in \mathcal{J}} \hat{X}_i\right) \\
&\leq 2(\kappa\beta_2)^{1/2} + \kappa\beta_2 \leq 2/3.
\end{aligned}
$$
(5.6)

Applying (3.23) to $\bar{W} - E\bar{W}$ yields

$$
\begin{aligned}
P(\bar{W} > z - 2) &\leq \frac{E|\bar{W} - E\bar{W} + 4|^4}{(z + 2 - E\bar{W})^4} \\
&\leq C(z+1)^{-4}\left(1 + \kappa^3 \sum_{i \in \mathcal{J}} E|X_i|^4 I(|X_i| \leq \tau)\right) \\
&\leq C(z+1)^{-4}\left(1 + \kappa^3 \tau \sum_{i \in \mathcal{J}} E|X_i|^3 I(|X_i| \leq \tau)\right) \\
&\leq C(z+1)^{-4}(1 + \kappa^2\beta_3).
\end{aligned}
$$
(5.7)

By (5.4)–(5.7) and the assumption that $\kappa\beta_2 \leq 1$, (5.3) is proved and hence the lemma. $\square$

LEMMA 5.2. *Assume $E|X_i|^p < \infty$ for some $2 < p \leq 3$. Then there exists an absolute constant $C$ such that, for $z \geq 0$,*

$$P(W > (z-1)/2) \leq C(1+z)^{-p}(1 + \kappa^{p-1}\gamma) \tag{5.8}$$

*and*

$$
\begin{aligned}
&\sum_{i \in \mathcal{J}} P(|Y_i| > (z+1)/4) + \kappa^3(1+z)^{-3}\beta_3 \\
&\quad + (1+z)^{-3}\beta_1 + \kappa^2\beta_2\{(1+z)^{-3} + P(W \geq (z-1)/2)\} \\
&\quad \leq C\kappa^p(1+z)^{-p}\gamma + C(1+z)^{-3}\kappa^{2p-1}\gamma^2,
\end{aligned}
$$
(5.9)

*where $\gamma = \sum_{i \in \mathcal{J}} E|X_i|^p$.*

PROOF. If $\kappa^{p-1}\gamma > (1+z)^{p-2}$, then (5.8) is trivial because $P(W > (z-1)/2) \leq 4(z+1)^{-2}E(W+1)^2 = 8(z+1)^{-2}$. To prove (5.8) for $\kappa^{p-1}\gamma \leq (1+z)^{p-2}$, let

$$
\begin{aligned}
\tilde{X}_i &= X_i I(|X_i| \leq (1+z)\tau) - EX_i I(|X_i| \leq (1+z)\tau), \\
\tilde{X}_i^* &= X_i I(|X_i| > (1+z)\tau) - EX_i I(|X_i| > (1+z)\tau), \\
\tilde{W} &= \sum_{i \in \mathcal{J}} \tilde{X}_i, \qquad \tilde{W}^* = \sum_{i \in \mathcal{J}} \tilde{X}_i^*.
\end{aligned}
$$

Then
$$P(W > (z-1)/2)$$
$$\leq P(\tilde{W} > (z-3)/4) + P(\tilde{W}^* > (z+1)/4)$$
$$\leq C(1+z)^{-4} E|\tilde{W}+1|^4 + 4(1+z)^{-1} E|\tilde{W}^*|$$
$$\leq C(1+z)^{-4}(1+E|\tilde{W}|^4) + C(1+z)^{-1}((1+z)\tau)^{-p+1}\gamma.$$

Similarly to (5.6), $\kappa^{p-1}\gamma \leq (1+z)^{p-2}$ implies $\mathrm{Var}(\tilde{W}) \leq 4$. Hence by (3.23),
$$E|\tilde{W}|^4 \leq C\left(1 + \kappa^3 \sum_{i \in \mathcal{J}} E|X_i|^4 I(|X_i| \leq (1+z)\tau)\right)$$
$$\leq C(1 + \kappa^{p-1}(1+z)^{4-p}\gamma).$$

By combining the above inequalities, (5.8) is proved.

We now prove (5.9). From (5.8), the Chebyshev inequality and the Hölder inequality, the left-hand side of (5.9) is bounded by
$$((z+1)/4)^{-p} \sum_{i \in \mathcal{J}} E|Y_i|^p + \kappa^3 (1+z)^{-3} \tau^{3-p} \gamma$$
$$+ (1+z)^{-3} \tau^{-p+2} \sum_{i \in \mathcal{J}} E|X_i|^{p-1}|Y_i| + C(1+z)^{-3}\kappa^2 \tau^{-p+2}\gamma(1+\kappa^{p-1}\gamma)$$
$$\leq C(1+z)^{-p}\kappa^p \gamma + C(1+z)^{-3}\kappa^{2p-1}\gamma^2.$$

This completes the proof of the lemma. □

PROOF OF THEOREM 2.5. Without loss of generality, assume $z \geq 0$. When $\kappa^{p-1}\gamma > 1$, (2.10) follows directly from (5.8). When $\kappa^{p-1}\gamma \leq 1$, then (2.10) is a consequence of (5.9) and the following inequality:

$$\begin{aligned}(5.10)\quad |F(z) - \Phi(z)| &\leq \sum_{i \in \mathcal{J}} P(|Y_i| > (z+1)/4) + C(1+z)^{-3}\beta_1 \\ &\quad + C\kappa^2 \beta_2 \{(1+z)^{-3} + P(W \geq (z-1)/2)\} \\ &\quad + C\kappa^2(1+z)^{-3}\beta_3.\end{aligned}$$

So it suffices to prove (5.10). It is clear, that for $z \geq 0$,
$$|P(W > z) - (1 - \Phi(z))| \leq P(W > (z-1)/2) + 16(1+z)^{-3}.$$

So (5.10) holds if $\kappa^2 \beta_2 > \frac{1}{28}$. It remains to consider the case that

$$(5.11) \quad \kappa^2 \beta_2 \leq \tfrac{1}{16}.$$

Let $\bar{W}$ be defined as in (5.1). By Lemma 5.1, it suffices to show that

$$(5.12)\quad \begin{aligned}|P(\bar{W} \leq z) - \Phi(z)| \\ \leq C(1+z)^{-3}\beta_1 + C\kappa(1+z)^{-3}\beta_2 + C\kappa^2(1+z)^{-3}\beta_3.\end{aligned}$$

Let $\sigma^2 = \text{Var}(\bar{W})$, $T = (\bar{W} - E\bar{W})/\sigma$. Then

(5.13)
$$\begin{aligned}|P(\bar{W} \leq z) - \Phi(z)| \\ = |P(T \leq (z - E\bar{W})/\sigma) - \Phi(z)| \\ \leq |P(T \leq (z - E\bar{W})/\sigma) - \Phi((z - E\bar{W})/\sigma)| \\ + |\Phi((z - E\bar{W})/\sigma) - \Phi(z)|.\end{aligned}$$

By the Chebyshev inequality,

(5.14) $$|E\bar{W}| \leq \tau^{-1} \sum_{i \in \mathcal{J}} EX_i^2 I(|X_i| > \tau) \leq \tfrac{1}{16}.$$

By (5.6), we have $1/3 < \sigma^2 < 2$, and moreover, similarly to (5.6),

(5.15)
$$\begin{aligned}|\sigma^2 - 1| &= 2\left|\sum_{i \in \mathcal{J}} E(\hat{X}_i W)\right| + \text{Var}\left(\sum_{i \in \mathcal{J}} \hat{X}_i\right) \\ &\leq 2\sum_{i \in \mathcal{J}} E|X_i Y_i| I(|X_i| > \tau) + \kappa \sum_{i \in \mathcal{J}} EX_i^2 I(|X_i| > \tau) \\ &= 2\beta_1 + \kappa\beta_2.\end{aligned}$$

Thus by (5.14) and (5.15),

(5.16)
$$\begin{aligned}|\Phi((z - E\bar{W})/\sigma) - \Phi(z)| \\ \leq |\Phi((z - E\bar{W})/\sigma) - \Phi(z/\sigma)| + |\Phi(z/\sigma) - \Phi(z)| \\ \leq C(1+z)^{-3}|E\bar{W}| + C(1+z)^{-3}|\sigma^2 - 1| \\ \leq C(1+z)^{-3}\kappa\beta_2 + C(1+z)^{-3}\beta_1.\end{aligned}$$

With $x = (z - E\bar{W})/\sigma$ $(> -1/2)$, we only need to show that

(5.17) $$|P(T \leq x) - \Phi(x)| \leq C\kappa^2 (1+x)^{-3}\beta_3.$$

Put

$$\xi_i = (\bar{X}_i - E(\bar{X}_i))/\sigma, \qquad \xi_{A_i} = \sum_{j \in A_i} \xi_j, \qquad \xi_{B_i} = \sum_{j \in B_i} \xi_j,$$

$$\hat{K}_i(t) = \xi_i(I(-\xi_{A_i} < t \leq 0) - I(0 < t < -\xi_{A_i})), \qquad K_i(t) = E\hat{K}_i(t).$$

By the definition of $\tau$, we have $|\xi_{A_i}| \leq \tfrac{1}{2}$ and $|\xi_{B_i}| \leq \tfrac{1}{2}$. If $(1+x)\kappa^2\beta_3 > 1$, then by (3.23),

$$\begin{aligned}|P(T > x) - (1 - \Phi(x))| \\ \leq (1+x)^{-4} E|1+T|^4 + (1+x)^{-4} \\ \leq C(1+x)^{-4}\left(1 + \kappa^3 \sum_{i \in \mathcal{J}} E|\xi_i|^4\right) \\ \leq C(1+x)^{-4}\left(1 + \kappa^3 \tau \sum_{i \in \mathcal{J}} E|\xi_i|^3\right) \\ \leq C(1+x)^{-4}(1 + \kappa^2\beta_3) \leq C(1+x)^{-3}\kappa^2\beta_3,\end{aligned}$$

which proves (5.17).

If $(1+x)\kappa^2\beta_3 \leq 1$, let $\alpha = 64\kappa^2\beta_3$. Also let $h_{x,\alpha}(w)$ be as in (4.4) and let $f(w) = f_{x,\alpha}(w)$ be the unique bounded solution of the Stein equation (4.3) with $x$ replacing $z$. Then by (4.1) and similarly to (4.2),

$$
\begin{aligned}
Ef'(T) &- ETf(T) \\
&= \sum_{i \in \mathcal{J}} E \int_{|t| \leq 1} (f'(T - \xi_{B_i} + t) - f'(T + t))\hat{K}_i(t)\,dt \\
&\quad + \sum_{i \in \mathcal{J}} E \int_{|t| \leq 1} (f'(T) - f'(T - \xi_{B_i} + t))K_i(t)\,dt \\
&:= R_1 + R_2.
\end{aligned}
\tag{5.18}
$$

Let $g(w) = (wf(w))'$ and let

$$R_{1,1} = \sum_{i \in \mathcal{J}} E\left|\int_{|t| \leq 1} \int_0^{-\xi_{B_i}} g(T+u)\,du\,\hat{K}_i(t)\,dt\right|,$$

$$R_{1,2} = \alpha^{-1}\sum_{i \in \mathcal{J}} E\left|\int_{|t| \leq 1} \int_0^{-\xi_{B_i}} I(x \leq T + u \leq x + \alpha)\,du\,\hat{K}_i(t)\,dt\right|.$$

Noting that

$$
\begin{aligned}
\left|f'_{x,\alpha}(w+t) - f'_{x,\alpha}(w+s) - \int_s^t g(w+u)\,du\right| \\
\leq \alpha^{-1}\left|\int_s^t I(x \leq w + u \leq x + \alpha)\,du\right|,
\end{aligned}
\tag{5.19}
$$

we have

$$|R_1| \leq R_{1,1} + R_{1,2}.$$

Let $\zeta_i = (\xi_i, \xi_{A_i}, \xi_{B_i})$. By Lemma 5.3,

$$
\begin{aligned}
R_{1,1} &\leq \sum_{i \in \mathcal{J}} E \int_{|t| \leq 1}\left|\int_0^{-\xi_{B_i}} E^{\zeta_i} g(T+u)\,du\right||\hat{K}_i(t)|\,dt \\
&\leq C(1+x)^{-3} \sum_{i \in \mathcal{J}} E \int_{|t| \leq 1} |\xi_{B_i}||\hat{K}_i(t)|\,dt \\
&\leq C(1+x)^{-3} \sum_{i \in \mathcal{J}} E|\xi_i \xi_{A_i} \xi_{B_i}| \\
&\leq C\kappa^2(1+x)^{-3}\beta_3,
\end{aligned}
$$

and by Proposition 3.3,

$$R_{1,2} \leq \alpha^{-1}\sum_{i \in \mathcal{J}} E \int_{|t| \leq 1}\left|\int_0^{-\xi_{B_i}} P^{\zeta_i}(x \leq T + u \leq x + \alpha)\,du\right||\hat{K}_i(t)|\,dt$$

$$\leq C\alpha^{-1}(1+x)^{-3}\sum_{i\in\mathcal{J}} E\int_{|t|\leq 1}(\kappa^2\beta_3+\alpha)|\xi_{B_i}||\hat{K}_i(t)|\,dt$$

$$\leq C\alpha^{-1}(1+x)^{-3}(\kappa^2\beta_3+\alpha)\sum_{i\in\mathcal{J}} E|\xi_i\xi_{A_i}\xi_{B_i}|$$

$$\leq C\kappa^2(1+x)^{-3}\beta_3.$$

This proves
$$|R_1|\leq C\kappa^2(1+x)^{-3}\beta_3.$$

Similarly we have
$$|R_2|\leq C\kappa^2(1+x)^{-3}\beta_3.$$

Hence we have

(5.20) $$|Eh_{x,\alpha}(T)-Nh_{x,\alpha}|\leq C\kappa^2(1+x)^{-3}\beta_3.$$

Finally, using the fact that $Eh_{x-\alpha,\alpha}(T)\leq P(T\leq x)\leq Eh_{x,\alpha}(T)$ and that $|\Phi(x+\alpha)-\Phi(x)|\leq C(1+x)^{-3}\alpha$ for $x>-1/2$, we have (5.17). This completes the proof of Theorem 2.5. □

It remains to prove Lemma 5.3 which was used above.

LEMMA 5.3. *Let* $\zeta=\zeta_i=(\xi_i,\xi_{A_i},\xi_{B_i})$ *and*

(5.21) $$g(w)=g_{x,\alpha}(w)=(wf_{x,\alpha}(w))'.$$

*Then for* $x>-1/2$ *and* $|u|\leq 4$, *we have*

(5.22) $$E^\zeta g(T+u)\leq C(1+x)^{-3}(1+\kappa^2\beta_3).$$

PROOF. From the definitions of $Nh_{x,\alpha}$, $f_{x,\alpha}$ and $g$, we have

$$Nh_{x,\alpha}=\Phi(x)+\frac{\alpha}{\sqrt{2\pi}}\int_0^1 se^{-(x+\alpha-\alpha s)^2/2}\,ds,$$

$$f_{x,\alpha}(w)=\begin{cases}\sqrt{2\pi}e^{w^2/2}\Phi(w)(1-Nh_{x,\alpha}), & w\leq x,\\ \sqrt{2\pi}e^{w^2/2}(1-\Phi(w))Nh_{x,\alpha}\\ \quad -\alpha e^{w^2/2}\int_0^{1+(z-w)/\alpha} se^{-(z+\alpha-\alpha s)^2/2}\,ds, & x<w\leq x+\alpha,\\ \sqrt{2\pi}e^{w^2/2}(1-\Phi(w))Nh_{x,\alpha}, & w>x+\alpha,\end{cases}$$

$$g(w)=\begin{cases}(\sqrt{2\pi}(1+w^2)e^{w^2/2}\Phi(w)+w)(1-Nh_{x,\alpha}), & w<x,\\ (\sqrt{2\pi}(1+w^2)e^{w^2/2}(1-\Phi(w))-w)Nh_{x,\alpha}+r_{x,\alpha}(w),\\ & x\leq w\leq x+\alpha,\\ (\sqrt{2\pi}(1+w^2)e^{w^2/2}(1-\Phi(w))-w)Nh_{x,\alpha}, & w>x+\alpha,\end{cases}$$

where

$$r_{x,\alpha}(w) = -we^{w^2/2} \int_w^{x+\alpha} \left(1 + \frac{x-s}{\alpha}\right) e^{-s^2/2} \, ds + \left(1 + \frac{x-w}{\alpha}\right).$$

For $x < w < x + \alpha$, we have, by integration by parts,

$$\begin{aligned} r_{x,\alpha} &\geq -e^{w^2/2} \int_w^{x+\alpha} s\left(1 + \frac{x-s}{\alpha}\right) e^{-s^2/2} \, ds + \left(1 + \frac{x-w}{\alpha}\right) \\ &= e^{w^2/2} \int_w^{x+\alpha} e^{-s^2/2} \, ds \geq 0. \end{aligned}$$

So $0 \leq r_{x,\alpha} \leq 1$ for $x < w < x + \alpha$. It can be verified that $0 < g(w) \leq C(1+|x|)$ for all $x$. Therefore, (5.22) holds when $-1 < x \leq 6$.

When $x > 6$, the proof of (5.22) is very similar to that of Lemma 5.2 of Chen and Shao (2001).

A direct calculation shows that, for $w \geq 0$,

(5.23) $$0 \leq \sqrt{2\pi}(1+w^2)e^{w^2/2}(1-\Phi(w)) - w \leq \frac{2}{1+w^3}.$$

This implies that $g \geq 0$, $g(w) \leq 2(1 - \Phi(x))$ for $w \leq 0$ and $g(w) \leq \frac{2}{1+w^3} + I(x < w < x+\alpha)$ for $w \geq x$; furthermore, $g$ is clearly increasing for $0 \leq w < x$. Therefore,

(5.24)
$$\begin{aligned} E^\zeta g(T+u) &= E^\zeta g(T+u) I(T+u \leq x-1) \\ &\quad + E^\zeta g(T+u) I(T+u \geq x) \\ &\quad + E^\zeta g(T+u) I(x-1 < T+u < x) \\ &\leq 2(1-\Phi(x)) + g(x-1) + 2(1+x^3)^{-1} \\ &\quad + P^\zeta(x \leq T+u \leq x+\alpha) \\ &\quad + E^\zeta g(T+u) I(x-1 < T+u < x) \\ &\leq C((1+x)^{-3} + x^2 e^{(x-1)^2/2}(1-\Phi(x))) + C(1+x)^{-3}\alpha \\ &\quad + E^\zeta g(T+u) I(x-1 < T+u < x) \\ &\leq C(1+x)^{-3}(1+\alpha) + E^\zeta g(T+u) I(x-1 < T+u < x). \end{aligned}$$

By Proposition 3.3 and the fact that $w - u \geq x - 5 > (x+1)/7$ for $w \geq x - 1 \geq 5$ and $|u| \leq 4$,

$$
\begin{aligned}
E^\zeta g(T+u)&I(x-1 < T+u < x) \\
&= \int_{x-1}^{x} -g(w)\, dP^\zeta(w < T+u < x) \\
&= g(x-1)P^\zeta(x-1 < T+u < x) \\
&\quad + \int_{x-1}^{x} g'(w) P^\zeta(w < T+u < x)\, dw \\
&\leq C(1+x)^{-3} + C(1+x)^{-3} \int_{x-1}^{x} g'(w)\{\alpha + (x-w)\}\, dw \\
&\leq C(1+x)^{-3} \left\{1 + \alpha g(x-) + \int_{x-1}^{x} (x-w)\, dg(w)\right\} \\
&\leq C(1+x)^{-3} \left\{1 + \alpha x + \int_{x-1}^{x} g(w)\, dw\right\} \\
&\leq C(1+x)^{-3} \{1 + x f_{x,\alpha}(x)\} \\
&\leq C(1+x)^{-3}.
\end{aligned}
\tag{5.25}
$$

Combining (5.24) and (5.25) yields (5.22). This completes the proof of Lemma 5.3. □

**Acknowledgment.** The authors thank the referee for valuable comments which led to substantial improvement in the presentation of the paper.

INSTITUTE FOR MATHEMATICAL SCIENCES  
NATIONAL UNIVERSITY OF SINGAPORE  
3 PRINCE GEORGE'S PARK  
SINGAPORE 118402  
SINGAPORE  
AND  
DEPARTMENT OF MATHEMATICS  
DEPARTMENT OF STATISTICS AND  
  APPLIED PROBABILITY  
NATIONAL UNIVERSITY OF SINGAPORE  
SINGAPORE 117543  
SINGAPORE  
E-MAIL: lhychen@ims.nus.edu.sg  

DEPARTMENT OF MATHEMATICS  
UNIVERSITY OF OREGON  
EUGENE, OREGON 97403  
USA  
AND  
DEPARTMENT OF MATHEMATICS  
DEPARTMENT OF STATISTICS AND  
  APPLIED PROBABILITY  
NATIONAL UNIVERSITY OF SINGAPORE  
SINGAPORE 117543  
SINGAPORE  
E-MAIL: qmshao@darkwing.uoregon.edu